\documentclass[12pt]{article}
\usepackage{anysize}
\marginsize{3cm}{1.5cm}{2cm}{2cm}
\usepackage{amsmath}
\usepackage{amsfonts}
\usepackage{amssymb}
\usepackage{diagrams}
\usepackage{theorem}
\DeclareMathAlphabet{\mathscr}{OT1}{pzc}{m}{it}
\usepackage{eucal}
\usepackage{verbatim}
\usepackage{graphicx}
\usepackage{indentfirst}
\usepackage{hyperref}
\usepackage[T1]{fontenc}
\normalfont

\makeatletter
\begingroup
\gdef\th@dotted{\normalfont\itshape
  \def\@begintheorem##1##2{%
        \item[\hskip\labelsep \theorem@headerfont ##1\ ##2.]}%
\def\@opargbegintheorem##1##2##3{%
   \item[\hskip\labelsep \theorem@headerfont ##1\ ##2\ (##3).]}}
\endgroup
\makeatother

\theoremstyle{dotted}

\newtheorem{thm}{Theorem}[section]
\newtheorem{lemma}[thm]{Lemma}
\newtheorem{prop}[thm]{Proposition}
\newtheorem{corr}[thm]{Corollary}

\newtheorem{introthm}{Theorem}

\makeatletter
\begingroup
\gdef\th@upshape{\normalfont
  \def\@begintheorem##1##2{%
        \item[\hskip\labelsep \theorem@headerfont ##1\ ##2.]}%
\def\@opargbegintheorem##1##2##3{%
   \item[\hskip\labelsep \theorem@headerfont ##1\ ##2\ (##3).]}}
\endgroup
\makeatother

\theoremstyle{upshape}

\newtheorem{opr}[thm]{Definition}
\newtheorem{rem}[thm]{Remark}
\newtheorem{ntn}[thm]{Notation}
\newtheorem{exm}[thm]{Example}
\newtheorem{const}[thm]{Construction}

\newtheorem{introopr}[introthm]{Definition}

\newcommand{\proof}[1][Proof. ]{\smallskip\noindent{\bf #1}}
\def\endproof{\hfill\ensuremath{\square}\par\medskip}

\makeatletter
\@addtoreset{equation}{section}
\makeatother

\newenvironment{introequation}{%

\begin{equation}}%
{\end{equation}}%

\newcommand{\sI}{{\mathscr I}}
\newcommand{\sJ}{{\mathscr J}}

\newcommand{\bG}{{\mathbb G}}
\newcommand{\bB}{{\mathbb B}}
\newcommand{\bC}{{\mathbb C}}
\newcommand{\bR}{{\mathbb R}}

\newcommand{\bD}{{\mathbb D}}

\newcommand{\bI}{{\mathbb I}}
\newcommand{\bJ}{{\mathbb J}}

\newcommand{\rB}{{\textrm B}}
\newcommand{\rR}{{\textrm R}}
\newcommand{\rL}{{\textrm L}}

\newcommand{\cB}{{\cal B}}

\newcommand{\cC}{{\cal C}}
\newcommand{\cD}{{\cal D}}

\newcommand{\cF}{{\cal F}}
\newcommand{\cN}{{\cal N}}
\newcommand{\cM}{{\cal M}}

\newcommand{\cW}{{\cal W}}

\newcommand{\cX}{{\cal X}}

\newcommand{\cE}{{\cal E}}
\newcommand{\cK}{{\cal K}}

\newcommand{\cO}{{\cal O}}

\renewcommand{\phi}{\varphi}

\newarrow {Into} {boldhook}--->
\newarrow {Mapsto} |--->
\newcommand{\Set}{{\rm \bf Set}}
\newcommand{\Fin}{{\rm \bf Fin}}
\newcommand{\SSet}{{\rm \bf SSet}}
\newcommand{\Cat}{{\rm \bf Cat }}
\newcommand{\Top}{{\rm \bf Top}}

\newcommand{\DGMod}{{\rm \bf DVect}}

\newcommand{\Ho}{{\rm Ho} \,}

\newcommand{\bfE}{\mathbf E}
\newcommand{\bF}{\mathbb F}
\newcommand{\bO}{\mathbb O}
\newcommand{\op}{{\sf op}}
\newcommand{\bd}{\mathbf d}
\newcommand{\bc}{\mathbf c}

\newcommand{\bi}{\mathbf i}
\newcommand{\bj}{\mathbf j}

\renewcommand{\lim}{\operatorname{\varprojlim}}
\newcommand{\colim}{\operatorname{\varinjlim}}

\newcommand{\Map}{\operatorname{Map}}
\newcommand{\Sect}{{\sf Sect}}
\newcommand{\PSect}{{\sf PSect}} 
\newcommand{\DSect}{{\sf DSect}}
 
\newcommand{\Cart}{{\sf Cart}} 
\newcommand{\ovr}{/ \! \! /}

\author{Edouard Balzin}
\title{Derived sections of Grothendieck fibrations and the problems of homotopical algebra}
\date{}

\begin{document}
\maketitle
\tableofcontents

\section*{Introduction}

\subsection*{$\mathbb E_n$-algebras}

The formalism of topological operads \cite{MAY} appeared as a way to describe the algebraic structure of $n$-fold loop spaces. A topological operad $\cO$ is a symmetric sequence of spaces $\{\cO(l)\}_{l \in \mathbb N}$, where each $\cO(l)$ should be thought as the space of operations with $n$ inputs and one output, supplied with composition maps $\cO(l) \times \cO(m) \to \cO(l+m-1)$ satisfying symmetry, associativity and unitality conditions. Of particular interest are little $n$-disk operads $\cD_n$ (defined as little cube operads in \cite{MAY}), for which $\cD_n(l)$ is (homotopy equivalent to) the configuration space of $l$ points in $\mathbb R^n$. Any $n$-fold loop space $X$ has a structure of an algebra over $\cD_n$, which means there are maps $\cD_n(l) \times X^l \to X$ satisfying certain conditions.

Denote by $\DGMod_k$ the category of chain complexes over a field $k$. Taking singular chains of $\cD_n(l)$ produces the  $\DGMod_k$-operad $\mathbb E_n$ (and indeed, $\cD_n$ are often called $\mathbb E_n$-operads in topological spaces). Algebras over such operads $\mathbb E_n$, that is, chain complexes $M$ together with action maps $\mathbb E_n(l) \otimes M^{\otimes l} \to M$, have been studied with great interest especially in the recent years \cite{LU}. An example of an $\mathbb E_2$-algebra is the cohomological Hochschild complex $CH(A)$ of a DG-algebra $A$, which appears in many settings, for example in two-dimensional topological conformal field theories \cite{COSTELLO}. It is important, however, to remark that $CH(A)$ is an $\mathbb E_2$ algebra \emph{up to quasi-isomorphism}: the Deligne Conjecture \cite{BER-FRE,MCCLSM,TAM} only implies that there exists an operad $\cO$ in $\DGMod_k$, quasi-isomorphic to $\mathbb E_2$, which acts on $CH(A)$. The proofs of this result involve, subsequently, a lot of combinatorial work to construct $\cO$, its action on $CH(A)$, and the (chain of) quasi-isomorphisms $\cO \cong \mathbb E_2$. 

The bulkiness of the formalism of operads comes from the fact that two quasi-isomorphic or homotopy equivalent operads can be of very different size and complexity yet describe equivalent structures. There is, however, a different approach to $\mathbb E_n$-algebras, and more generally to structures related to configuration spaces, which relies on the machinery of factorisation algebras of \cite{BD}. A factorisation algebra $\cF$ over a space $X$ consists of, roughly speaking, a presheaf $\cF_n$ of complexes of vector spaces on $X^n$ for each power $n$, together with additional structure. First, there is a map
\begin{introequation}
\label{introfactorisationdiag}
\Delta^n_* \cF_n \longrightarrow \cF_1
\end{introequation}%
between the restriction $\Delta_n^* \cF_n$ of $\cF_n$ along the smallest diagonal $\Delta_n: X \to X^n$, and $\cF_1$. Second, if we denote by $i_n:U_n \subset X_n$ the complement $\{(x_i) | x_k \neq x_l \}$ to all diagonals, then there are factorisation maps 
\begin{introequation}
\label{introfactorisationmaps}
i_n^* \cF_n \longrightarrow \cF_1 \boxtimes ... \boxtimes \cF_1 
\end{introequation}%
between the restriction of $\cF_n$ to $U_n$ and the external product of $\cF_1$ \cite{BD}, which are required to be quasi-isomorphisms. When $X$ is a $k$-disk, one can prove that $\mathbb E_n$-algebras correspond to those factorisation algebras on $X$ which are moreover constructible (which means that $\cF_n$ is locally constant on the strata for the standard stratification of $X^n$). 

The notion of factorisation algebra has proved its use and is arguably more natural and canonical than algebras over topological operads. In dimension $2$, factorisation algebras are particularly useful, as constructible sheaves on a two-disk $D$ can be shown to depend on the data which can be described using a version of the fundamental groupoid of the configuration space $D^n$. Comparing this to operads, one may thus wonder if there is a general ``homotopic-algebraic'' formalism which does not suffer from the noncanonicity issues of topological operads, and naturally reproduces factorisation algebra approach to a variety of algebraic structures.

\subsection*{The approach of Segal}
For the context of loop spaces, an approach alternative to operads does exist and is very useful in practical applications. In \cite{SEG}, Graeme Segal introduced the notion of a $\Gamma$-space. Denote by $\Fin_*$ the category of finite sets and partially defined maps: a map $S \to T$ in $\Fin_*$ is a map of sets $U \to T$ defined on a subset $U \subset S$. A $\Gamma$-space $A$ is then defined as a functor
\begin{diagram}[small]
\Fin_* & \rTo^A & \Top
\end{diagram}
to the category of topological spaces $\Top$, satisfying Segal conditions. To formulate them, fix a one-element set $1$. For a set $S$ and an element $x \in S$, we have the corresponding partial map $\rho_x: S \to 1$ defined on the subset $\{ x \}$. The Segal conditions say that for each $S \in \Fin_*$, the induced map
\begin{introequation}
\label{introsegalconditions}
\begin{diagram}[small]
A(S) & \rTo^{\prod_{x \in S} A(\rho_x)} &   A(1)^S
\end{diagram}
\end{introequation}%
is a homotopy equivalence of topological spaces.

For each $S \in \Fin_*$ there is one map $\pi_S: S \to 1$ defined on the whole of $S$. We can consider the following span
\begin{introequation}
\label{introspan}
\begin{diagram}[small,nohug]
		&							& A(S)	&					&	\\
		&	\ldTo^{\prod A(\rho_x)}	&		& \rdTo^{A(\pi_S)}	&		\\
 A(1)^S 	&							&		&					&	A(1).
\end{diagram}
\end{introequation}%
The Segal conditions imply that this span represents a morphism in the homotopy category $\Ho \Top$. A choice of a homotopy inverse for the left map gives, non-canonically, a multiplication operation $m_S: A(1)^S   \to A(1)$. After projecting to $\Ho \Top$ one can check that $A(1)$ becomes a commutative monoid. 

However, a $\Gamma$-space $A$ carries more information than the homotopy monoid $A(1)$. Segal, just like May with topological operads, used $\Gamma$-spaces to describe infinite loop spaces and his delooping machinery. From the modern perspective, a $\Gamma$-space is a proper description of a homotopy coherent commutative monoid in topological spaces. In particular, $\Gamma$-spaces describe the same structure as $\mathbb E_\infty$ \cite{MAY} algebras in $\Top$.

Instead of $\Fin_*$ we can consider the opposite of the usual category of simplices $\Delta$. A suitable modification of definitions then permits to model homotopy coherent monoids with no commutativity as certain simplicial spaces $\Delta^\op \to \Top$. Explicit examples include ordinary loop spaces. Moreover, the work of \cite{BAR} implies that there are categories $C$ such that $n$-fold loop spaces --- examples of $\cD_n$-algebras --- can be also modelled as Segal-type objects $C \to \Top$ for a suitable choice of the category $C$. In place of $\Top$, one can consider any homotopical category, that is a category $\cM$ with a subcategory of weak equivalences $\cW$, such that $\cM$ has (homotopy) products, and define Segal objects as functors to $\cM$ with maps (\ref{introsegalconditions}) being weak equivalences.         

The Segal space approach contrasts with operadic approach in that multiplicative operations $m_S: A(1)^S \to (1)$ for a $\Gamma$-space $A$ are not defined canonically and instead are constructed using the properties of $A$, while specifying a model $\cO$ for $\mathbb E_\infty$-operad in $\Top$ and an algebra over it means supplementing a lot of structure. In particular, for a $|S|$-element set $S$, $A(S)$ need not to be equal to $\cO (|S|) \times A(1)^S$. The information about multiplication properties in Segal formalism is thus entirely contained in the category $\Fin_*$. There is much less arbitrary choice left available, and one might hope it would be easier to construct and work with Segal structures rather than with operadic structures. Moreover, there is a great similarity between Segal $\Gamma$-spaces and factorisation algebras: for a factorisation algebra $\cF$, the maps (\ref{introfactorisationdiag}) and (\ref{introfactorisationmaps}) provide, after passing to stalks, spans just like (\ref{introspan}). 

However, if we attempt to extend the formalism of Segal objects to chain complexes, we immediately run into difficulties. To produce maps like (\ref{introsegalconditions}) in the $\Gamma$-space picture we used the universal property of Cartesian product $\times$ which is not satisfied by the tensor product $\otimes_k$ of $\DGMod_k$.

To deal with this issue, one observes \cite{SEG,TV} that any symmetric monoidal category $\cM$ is a weakly commutative monoid object in $\Cat$, the category of categories of suitable size. It is then, up to an equivalence, described by a `$\Gamma$-category' $M: \Fin_* \to \Cat$, with maps (\ref{introsegalconditions}) being equivalences of categories. In order not to choose an equivalence between $\cM$ and $M(1)$, one has to work either with pseudofunctors from $\Fin_*$ to $\Cat$, or equivalently, with Grothendieck opfibrations \cite{SGA1,VIST} (also called coCartesian fibrations) over $\Fin_*$: either notion encodes a weak $\Fin_*$-indexed family of categories.  

The way \cite{LU} to directly produce a Grothendieck opfibration out of a symmetric monoidal category $\cM$ with monoidal product $\otimes$ is as follows. Define $\cM^\otimes$ to be the category with objects $(S, \{M_x\}_{x \in S})$ where $S \in \Fin_*$ and each $M_x$ is an object of $\cM$. A morphism $(S, \{M_x\}_{x \in S}) \to (T, \{N_y\}_{y \in T})$ consists of a partially defined map $f:S \to T$, and for each $y \in T$, of a morphism $\otimes_{x \in f^{-1}(y)} M_x \to N_y$; when $f^{-1}(y)$ is empty, the monoidal product over it is the unit object. The compositions can then be defined with the help of the coherence isomorphisms for the product $\otimes: \cM \times \cM \to \cM$ and the unit object. The natural functor $p: \cM^\otimes \to \Fin_*$ is a Grothendieck opfibration, which means that the assignment $S \mapsto p^{-1}(S) = \cM^S$ is functorial in a weak but coherent way. 

Given a monoid object $A \in \cM$, we may define a section $\Fin_* \to \cM^\otimes$ of $p: \cM^\otimes \to \Fin_*$ as $S \mapsto (S,\{P_x\})$ with each $P_x = A$. Sections of this type can be characterised by putting suitable normalisation conditions on sections of $p$. However, there is still no evident way to write diagrams for Segal conditions using the language of sections. If we take a map $f: S \to T$ in $\Fin_*$, the value of a section $A$ on $f$ is determined by a map $f_! A(S) \to A(T)$ in $\cM^T$, where $f_!: \cM^S  \to \cM^T $ is the functor 
\begin{introequation}
\label{introtransitionfunctor}
f_!: (S,\{P_x\}_{x \in S}) \mapsto (T, \{Q_y\}_{y \in T}), \, \, \, Q_y = \otimes_{x \in f^{-1} y} P_x.
\end{introequation}%
In particular, for the map $\pi_S:S \to 1$ with $\pi_S^{-1}(1) = S$, we directly get multiplication operations $A(1)^{ \otimes S} = (\pi_S)_! A(S) \to A(1)$. To have a Segal description, we would instead like to have an object which produces diagrams like the following:
\begin{introequation}
\label{introsegcondtensor}
\begin{diagram}[small,nohug]
		&							& A_{\pi_S} 	&					&	\\
		&	\ldTo	&		& \rdTo	&		\\
A^{\otimes S} = (\pi_S)_! A(S) 	&							&		&					&	A(1).
\end{diagram}
\end{introequation}%
The left map may then be required to be a weak equivalence if $\cM$ has such. 

More generally, in place of $\cM^\otimes \to \Fin_*$ we can consider general Grothendieck opfibrations $\cE \to \cC$ and attempt to study their sections, in a homotopy sense.  While it is true that higher-categorical approaches to this problem exist \cite{LU,LUHTT}, the resulting machinery is extremely unwieldy, and one may wonder if the problem can be tackled using the means of classical categorical homotopy theory. We thus ask if there is a way to reproduce Segal approach, defining objects which, evaluated on a map $f:c \to c'$ in $\cC$, give spans of the form $$f_! A(c) \longleftarrow A_f \longrightarrow A(c')$$ (here $f_!: \cE(c) \to \cE(c')$ is the transition functor induced by the opfibration property of $\cE \to \cC$). Such objects would serve as a homotopical model for the sections of the opfibration $\cE \to \cC$.

\subsection*{Derived sections}

The purpose of this paper is to lay down a basic formalism which would answer the question posed in the last paragraph. In this paper we define \emph{derived} sections of opfibrations with weak equivalences, which in particular produce diagrams like (\ref{introsegcondtensor}) above, and prove a technical result which permits constructing examples of derived sections. 

Let us briefly describe our formalism. Recall first that for a category $\cC$, its simplicial replacement (a terminology choice inspired by \cite{BK}) $\bC$ is defined as the category whose objects are composable sequences $\bc_{[n]}=c_0 \to ... \to c_n$ of arrows of $\cC$ of arbitrary finite length $n \geq 0$. A morphism between $\bc_{[n]}=c_0 \to ... \to c_n$ and $\bc'_{[m]}=c'_0 \to ... \to c'_m$ consists of an order-preserving map of ordinals $a: [m] \to [n]$ (here $[i]$ denotes a totally ordered set of $i+1$ elements $0,1,...,i$) such that $c_{a(k)} = c'_k$ for $0 \leq k \leq m$. If, as before, one denotes by $\Delta$ the standard category of simplices, then $\bC$ is the opposite of the simplex category of the nerve $N \cC : \Delta^\op \to \Set$. The assignments $(c_0 \to ... \to c_n) \mapsto c_0$ or $(c_0 \to ... \to c_n) \mapsto c_n$ determine two `head' and `tail' functors $h: \bC \to \cC$ and $t: \bC \to \cC^\op$.   

Consider a functor $F: \bC \to \cM$ where $\cM$ is an arbitrary category with weak equivalences $\cW$. If we take a morphism $f: c \to c'$ of $\cC$, we then can consider the following span in $\bC$:
\begin{introequation}
\label{introspanbase}
\begin{diagram}[small,nohug]
		&							& c \stackrel f \to c' 	&					&	\\
		&	\ldTo	&		& \rdTo	&		\\
c 	&							&		&					&	c'.
\end{diagram}
\end{introequation}%
Evaluating $F$ on this diagram gives the corresponding span in $\cM$:
\begin{introequation}
\label{introspanfunctor}
\begin{diagram}[small,nohug]
		&							& F(c \stackrel f \to c') 	&					&	\\
		&	\ldTo	&		& \rdTo	&		\\
F(c) 	&							&		&					&	F(c').
\end{diagram}
\end{introequation}%
If one requires that $F(c) \longleftarrow F(c \stackrel f \to c')$ is an isomorphism, then the span (\ref{introspanfunctor}) defines a map from $F(c)$ to $F(c')$, which we denote as $F(f)$.
It then makes sense to ask if $F(gf) = F(g) F(f)$ for a composable pair of arrows $c \stackrel f \to c' \stackrel g \to c''$, or whether $F(id_c) = id_{F(c)}$. Both those conditions will be satisfied if $F$ sends to isomorphisms those maps of $\bC$ which have the form
$$
(c_0 \to ... \to c_k \to ... \to c_n) \longrightarrow (c_0 \to ... \to c_k)
$$
for $0 \leq k \leq n$ (that is, those maps which are determined by the inclusion of $[k]$ as first $k+1$ elements of $[m]$). Such maps of $\bC$ are called \emph{anchor maps} in this article. We observe that such a functor $F$ factors uniquely as $\bar F \circ h$, where $\bar F: \cC \to \cM$ is a functor from the original category $\cC$. 

If $F$ sends the anchor maps of $\bC$ to weak equivalences of $\cM$,  spans like (\ref{introspanfunctor}) define morphisms in $\Ho \cM$, the localisation of $\cM$ with respect to its weak equivalences. We may view such a functor $F: \bC \to \cM$ as a weakening of the notion of a functor from $\cC$  to $\cM$, with spans obtained from objects $c_0 \to ... \to c_n$ of greater length ensuring the coherence of compositions.

Assume now given an opfibration $p: \cE \to \cC$. Moreover, we assume that each fibre $\cE(c):= p^{-1}(c)$ has weak equivalences, and for each map $f: c \to c'$, the functor $f_!: \cE(c) \to \cE(c')$ induced by the opfibration property, preserves those weak equivalences. Then there exists a functor $p_\bC: \bfE \to \bC$, such that $\bfE(\bc_{[n]}) := p_\bC^{-1} (\bc_{[n]}) \cong \cE(c_n)$, and that for each $\alpha: \bc_{[n]} \to \bc'_{[m]}$, there is a naturally induced functor\footnote{Unlike $p$, the functor $p_\bC$ is a Grothendieck \emph{fibration} and describes a contravariant family over $\bC$.} $\bfE(\bc'_{[m]}) \to \bfE(\bc_{[n]})$ isomorphic to $t(\alpha)_!: \cE(c'_m) \to \cE(c_n)$.  

\begin{introopr}[Definition \ref{psectdef}]
	A \emph{presection} $B$ of $p: \cE \to \cC$ is a section $B: \bC \to \bfE$ of the functor $p_\bC$.
\end{introopr}
Presections form a category $\PSect(\cC,\cE)$ which can be equipped with weak equivalences if $p$ has such. A presection $B$ acting on spans like (\ref{introspanbase}) produces this diagram in $\cE(c')$:
\begin{diagram}[small,nohug]
		&							& B(c \stackrel f \to c') 	&					&	\\
		&	\ldTo	&		& \rdTo	&		\\
f_! B(c) 	&							&		&					&	B(c').
\end{diagram}
If the left map of this span and all other produced by applying $B$ to the anchor maps of $\bC$, are isomorphisms, then one can prove that $B$ defines an ordinary section $\cC \to \cE$ of the original opfibration $p: \cE \to \cC$. 

\begin{introopr}[Definition \ref{dsectdef}]
	A presection $B$ is called a \emph{derived section}, if $B$ takes anchor maps to weak equivalences (precisely to weakly Cartesian maps, see Definition \ref{wcart}).
\end{introopr}

The derived sections of $\cE \to \cC$ form a subcategory $\DSect(\cC,\cE) \subset \PSect(\cC,\cE)$ with induced weak equivalences. Denote by $\Ho \DSect (\cC,\cE)$ its localisation along the weak equivalences.

\subsection*{$\Delta$-categories and the direct image result}

The standard way to work with localisations like $\Ho \DSect(\cC,\cE)$ is to use model categories. However, in examples such as $\DGMod_k^\otimes \to \Fin_*$, while the fibres $\DGMod_k^\otimes(S) = \DGMod_k^S$ are model categories, the transition functors (\ref{introtransitionfunctor}) $f_!:\DGMod_k^S \to \DGMod_k^T$ do not usually preserve limits or colimits: in the basic case when $f: \{x,y\} \to 1$ is the map sending $x$ and $y$ to the element of $1$, $f_!:\DGMod_k^{\{x,y\}} = \DGMod_k \times \DGMod_k \to \DGMod_k$ is the tensor product $\otimes$, which preserves neither products nor direct sums\footnote{It is true that $\otimes$ preserves directed colimits, but this turns out to be insufficient in practice.}. This makes impossible applying the existing techniques \cite{H-S} for opfibrations and putting a model structure on $\PSect(\cC,\cE)$, let alone derived sections. 

As a consequence, for example, it is difficult to understand the behaviour of $\DSect(\cC,\cE)$ under the base change. Namely, given a functor $F: \cD \to \cC$, denote by $F^* p: F^* \cE \to \cD$ the pull-back of the opfibration $p: \cE \to \cC$. The functor $F^* p$ is again an opfibration, with fibrewise weak equivalences. There is an induced functor $\bF^* : \PSect(\cC,\cE) \to \PSect(\cD,F^* \cE)$, it preserves weak equivalences and respects derived sections. For a given $F$, we may ask if $\bF^*$ admits a homotopy adjoint, or is homotopically full and faithful when restricted to derived sections. And unfortunately, answering these questions without any additional structure on $p$ seems impossible. 


To deal with this problem, we propose the approach of homotopical $\Delta$-categories, which axiomatises the notion of geometric realisation of simplicial objects. Admittedly, this is a fairly restrictive technical tool specifically adapted to our setting, and the notion of derived section does not depend on $\Delta$-categories and can be studied using other means. Nonetheless, $\Delta$-categories permit a variety of constructions involving realisations of simplicial objects, which, together with homotopy coproducts, are known to produce all homotopy colimits in other homotopical settings \cite{LUHTT}. Examples of $\Delta$-categories include $\DGMod_k$, simplicial vector spaces $\Delta^\op \mathbf  {Vect}_k$ and some other categories which admit a sufficiently reasonable action of simplicial sets. For $\DGMod_k$ and $\Delta^\op \mathbf  {Vect}_k$, the $\Delta$-category structure interacts with the tensor product, giving rise to what we call homotopical $\Delta$-opfibrations. 

With the use of this technology, we are able to construct, for each $F: \cD \to \cC$ and a homotopical $\Delta$-opfibration $\cE \to \cC$, a weak equivalence preserving functor $$\bF_!: \PSect(\cD,F^* \cE) \to \PSect(\cC,\cE),$$ defined as an elaborate version of the bar construction \cite{RIEHL}. This functor is not adjoint to $\bF^*$ and does not, in general, preserve derived sections. However, $\bF_!$ comes together with additional structure, for example, with a natural transformation $\epsilon :\bF_! \bF^* \to id$ on the level of localisation $\Ho \DSect(\cC,\cE)$, which behaves like the counit of an adjunction data, see Proposition \ref{propositionadjunction} for details. These data are then used to prove some extra properties for the pair $\bF_!,\bF^*$ under some additional assumptions on $F$, which we explain below. 

While our construction of the functor $\bF_!$ uses $\Delta$-categories, admittedly, most of it depends only on the combinatorics of simplicial replacements. Thus if one finds another way to compute homotopy colimits in the setting of fibrations (say, with the use of model or higher categories), our construction can then be largely carried over. We would also like to note that our approach is similar in spirit to that of Costello \cite{COSTELLO}, who constructs a derived equivalence by providing explicitly two functors together with natural maps which become isomorphisms on the level of localisation. However, Costello's construction is rather ad hoc; we attempt to be more systematic.

\subsection*{Resolutions in the context of derived sections}

We then restrict our attention to one specific class of examples which are motivated by algebraic geometry \cite{KUZNETSOV, LUNTS} and called resolutions in this paper. A functor $F: \cD \to \cC$ is a \emph{resolution} if it is an opfibration and for each $c \in \cC$ the fibre $\cD(c) := F^{-1} (c)$ has contractible nerve. For an example of a resolution, consider a finite CW-complex $Y$ of homotopy type $K(G,1)$ and denote by $\rB G$ the fundamental groupoid of $Y$. Take $I$ to be the partially ordered set associated to a chosen regular cellular decomposition of $Y$. Choosing a point (say, the centre) of each cell of $I$ and connecting these points by paths when one cell is included in the other defines a functor 
\begin{introequation}
\label{introfresolut}
F:I \to \rB G
\end{introequation}%
which is (equivalent to) a resolution: it is an opfibration (up to an equivalence) since the category $\rB G$ is a groupoid, and its fibres are equivalent to cellular decompositions of the universal cover of $Y$, which is contractible by assumption. The functor $F$ induces the pull-back functor $F^*: \cD(\rB G, k) \to \cD(I,k)$, where $\cD(\rB G, k)$ and $\cD(I,k)$ are the derived categories of functors from $\rB G$ and $I$ to $\DGMod_k$ correspondingly (note that $\cD(\rB G, k)$ is the same thing as $Loc(Y,k)$, the derived category of complexes of locally constant sheaves on $Y$). One can prove that $F^*$ is full and faithful, with its image consisting of those functors $I \to \DGMod_k$ which are 'locally constant', in the sense that they send all morphisms of $I$ to quasi-isomorphisms. We also see that $\cD(I,k)$ is a good object: it is the category of modules over the finite-dimensional algebra generated by $I$. In particular, it is simple to construct objects of $\cD(I,k)$, which, if locally constant, can provide examples of $G$-representations. 


An example of a functor (\ref{introfresolut}) arises for the configuration spaces $\cD_2^n$ of $n$ points on a $2$-disk, which are the classifying spaces for $n$-braid groups $Br_n$, and which admit interesting cellular decompositions using the combinatorics of planar trees \cite{KON-SOB}. This example is of particular importance for Deligne conjecture and factorisation algebras.

For resolutions, the result we prove in this paper is the following. 

\begin{introthm}[Theorem \ref{maintheorem1}]
Given a homotopical $\Delta$-opfibration $\cE \to \cC$ and a resolution $F: \cD \to \cC$, the induced pull-back functor for derived sections, $\bF^* : \DSect(\cC,\cE) \to \DSect(\cD,F^* \cE)$, is homotopically full and faithful.
\end{introthm}

 In some specific cases we may also characterise (Theorem \ref{maintheorem2}) the essential homotopical image of $\bF^*$. It consists of those derived sections which are locally constant (in the weak sense) when restricted to the simplicial replacements $\bD(c)$ of fibres $\cD(c)$. In the case when our opfibration $\cE \to \cC$ is the constant opfibration $\DGMod_k \times \cC \to \cC$, we reproduce the result for derived categories discussed above. In general, resolutions serve as a testing case, indicating that $\DSect$ is a reasonable thing to consider. For the example of $\DGMod_k^\otimes \to \Fin_*$, resolutions can be used to provide the proof of Deligne conjecture without operadic considerations; the details were investigated by the author in his thesis and will be published in the near future. 

\subsection*{Organisation of the paper} In the first section, we introduce the formalism of homotopical $\Delta$-categories which we need in order to treat the issues arising from the breakdown of model-categorical formalism. The content of this section is not genuinely new and is somehow present in the folklore. For instance, the geometric realisation functor for homotopical categories, in the setting which goes beyond simplicial model categories, has been considered in \cite[Appendix]{BER-MOER}.

In the second section, we recall some of the basic notions and constructions related to the theory of Grothendieck (op)fibrations, including the less known constructions of transpose and power (op)fibrations. Since Grothendieck opfibrations are natural tools for encoding the notion of families of categories, we introduce a class of suitably structured opfibrations, called homotopical $\Delta$-opfibrations, which formalise the notion of a covariant family of $\Delta$-categories.

In the third section, we introduce simplicial replacements and then use them to define derived sections.  The fourth section deals with the construction of the pushforward functor $\bF_! $ and the map $\epsilon: \bF_! \bF^* \to id$, the data which one can use to verify if the `right adjoint' $\bF^*$  is full and faithful. 

Finally, the fifth section consists of the analysis of the case of a resolution, stating the Theorems \ref{maintheorem1} and \ref{maintheorem2} and outlining their proof. It is proven that in this case, the inverse image on the derived sections is full and faithful on the homotopy level. In addition, under mild assumptions we can characterise the essential image of the inverse image functor $\bF^*$.

\subsection*{Acknowledgements} The author is enormously grateful to his co-advisors, Dmitry Kaledin and Carlos Simpson, for their immense support and patience. In particular, Dmitry's insights on the problem were and are of great importance for the project. In the course of research, the author benefited a lot from conversations with Michael Batanin, Clemens Berger, Emily Riehl, Bruno Vallette and Gabriele Vezzosi. A separate thanks to each of the referees for providing comments of enormous value, which, hopefully, led to the improvement of this text. This work has been carried out in and with the help of University of Nice and Higher School of Economics, and the author enjoyed the friendly atmosphere of quite a few conferences that took place in Paris, Nice, Copenhagen, Lausanne and Yaroslavl. The author is grateful for the aid of TOFIGROU and HOGT ANR grants, which helped to cover some of his travel expenses. This paper uses Paul Taylor's diagrams package.

\section{Generalities on geometric realisation}

Apart from fixing notation, this section sets up the formalism of homotopical $\Delta$-categories. A $\Delta$-category $\cM$ is a category together with a multiplicative action of the simplex category $\Delta$ (see Definition \ref{deltacat}), and for homotopical $\Delta$-categories this action respects weak equivalences both in $\cM$ and simplicial sets, and provides us with a convenient version of geometric realisation functor. 

Our primary motivation for introducing such formalism is the structure on the category of chain complexes of vector spaces $\DGMod_k$. It is not a simplicial monoidal model category, but nonetheless comes with a $\Delta$-action given by tensor products with the chain complex of $n$-simplexes and an associated version of geometric realisation functor (see \cite{BER-MOER} for an alternative discussion of this example). 

\subsection{Homotopical categories}

\begin{ntn}
\label{notationdelta}
For any category $\cC$, $x \in \cC$ means that $x$ is an object of $\cC$. We also write $f \in \cC$ for morphisms $f: x \to y$ of $\cC$ if there is no confusion. The set of morphisms between two objects $x,y$ of $\cC$ is denoted $\cC(x,y)$. The category of functors $Fun(I,\cM)$ between two categories $I$ and $\cM$ is often denoted as $\cM^I$. Sometimes, given an object $x \in \cC$, we denote again by $x$ the functor from the terminal category to $\cC$ which picks out $x$.

From now on, $\Delta$ denotes the usual category of simplexes, i.e. the full subcategory of the category of small categories $\Cat$ spanned, for $n \geq 0$, by categories $[n]$ with $n+1$ objects $0,...,n$ and exactly one morphism from $i$ to $j$ whenever $i \leq j$.

By $\SSet= Fun(\Delta^\op, \Set)$ we denote the category of simplicial sets. We often identify $\Delta$ with its image in $\SSet$ by the Yoneda embedding
\begin{equation}
\label{deltayoneda}
\Delta^\bullet: \Delta \to \SSet = Fun(\Delta^\op, \Set), \, \, \, [n] \mapsto \Delta^n := \Delta^\bullet([n]) = \Delta(-,[n]).
\end{equation}

For a simplicial object $X: \Delta^\op \to \cM$ in a category $\cM$, denote $X_n:= X([n])$ for any $[n] \in \Delta$, and similarly, for bisimplicial objects $Y : \Delta^\op \times \Delta^\op \to \cM$, we write $Y_{nm} = Y(([n],[m]))$. We also write $\Delta^\op \cM:=Fun(\Delta^\op, \cM)$, $(\Delta \times \Delta^\op)\cM :=Fun(\Delta \times \Delta^\op, \cM)$ and so on.
\end{ntn}

\begin{opr}
A \emph{homotopical category} is a pair $(\cM,\cW)$ where $\cM$ is a category and $\cW$ is a subcategory of $\cM$ which contains all objects and isomorphisms.
We moreover require that for a composable pair of morphisms $f, g$ of $\cM$, if any two elements of $\{f,g, gf\}$ are in $\cW$, then the third one is in $\cW$ as well. 

We call $\cW$ the category of \emph{weak equivalences}. A morphism $f: x \to y$ of $\cM$ is a \emph{weak equivalence} if it belongs to $\cW$.
\end{opr}

\begin{opr}
\label{localisationdefinition}
For a homotopical category $(\cM, \cW)$ its \emph{localisation} \cite{DHKS,HOVEY} $\cW^{-1} \cM$ is the category together with a functor $p: \cM \to \cW^{-1} \cM$ such that any functor $F: \cM \to \cN$ which sends maps of $\cW$ to isomorphisms of $\cN$, factors through $p$ up to a canonical isomorphism. 
\end{opr}

We also denote $\cW^{-1} \cM$ by $\Ho \cM$. 

The existence of localisation for homotopical categories $(\cM,\cW)$ where $\cM$ is not small is a known set-theoretical issue, and one can check that for the examples of interest, in this paper it will not arise.

\begin{opr}
Given two homotopical categories $(\cM,\cW_\cM)$, $(\cN,\cW_\cN)$ a functor $F: \cM \to \cN$ is \emph{homotopical} iff $F(\cW_M) \subset \cW_N$. Equivalently, $F$ takes weak equivalences of $\cM$ to weak equivalences of $\cN$.
\end{opr}

Any homotopical functor $F: \cM \to \cN$ produces to a functor $\overline F: \Ho \cM \to \Ho \cN$.

\begin{exm}
Some well known examples of homotopical categories are
\begin{itemize}
\item the category $\SSet$ of simplicial sets which can be equipped with a homotopical structure by defining $\cW$ to be the subcategory of weak homotopy equivalences \cite{G-J} of simplicial sets,
\item the category $\DGMod_k$ of unbounded chain complexes over a field $k$, with $\cW$ being the subcategory of quasi-isomorphisms \cite{HOVEY},
\item any model category $\cM$ with its subcategory of weak equivalences $\cW$.
\end{itemize}
\end{exm}

\begin{opr}
\label{twooutofsix}
A subcategory $\cW \subset \cM$ satisfies the \emph{two-out-of-six} property, if given three maps in $\cM$ denoted $f, g, h$, so that they are composable with compositions $gf, hg, hgf$, if $gf, hg$ are in $\cW$, then all other maps $f, g, h$ and $hgf$ are in $\cW$.
\end{opr}

The subcategory of isomorphisms in any category satisfies two-out-of-six. The subcategory of weak equivalences in any model category satisfies two-out-of-six as well \cite{DHKS}.
 
\begin{opr}
For $I \in \Cat$ and a homotopical category $(\cM,\cW)$, the standard homotopical structure $(\cM^I,\cW_I)$ on the category of functors $A:I \to \cM$ consists of those natural transformations $\alpha: A \to B$ which are valued in the maps of $\cW$. That is, for each $i \in I$, the map $\alpha(i): A(i) \to B(i)$ is a weak equivalence. \end{opr}

\subsection{Tensors and $\Delta$-categories}
 
 The purpose of this subsection is to outline our formalism which has a built-in version of geometric realisation, and works well with the Grothendieck fibrations introduced later in the text. Denote by $\delta: \Delta \to \Delta \times \Delta$ the diagonal functor for $\Delta$.
 
\begin{opr}
\label{deltacat}
A $\Delta$-structure on a category $\cM$ consists of 
\begin{enumerate}
    \item a functor 
    $$
    \otimes: \Delta \times \cM \to \cM, \, \, \, ([n],x) \mapsto \Delta^n \otimes x,
    $$
    \item a natural transformation $diag$ depicted as a 2-square
    \begin{diagram}[small]
    \Delta \times \cM & \rTo^\otimes & \cM \\
    \dTo<{\delta \times id} & 	 \Downarrow diag 	& \uTo>\otimes \\
    \Delta \times \Delta \times \cM & \rTo_{id \times \otimes} & \Delta \times \cM \\   
    \end{diagram}
        \item a natural isomorphism of $\cM$-endofunctors: $\Delta^0 \otimes - \stackrel{\sim}{\to} id_\cM$.
\end{enumerate}

These data should satisfy the obvious coassociativity and counitality identities. A category $\cM$ with a $\Delta$-structure is called a \emph{$\Delta$-category} if $\cM$ is cocomplete and the functor $\otimes$ preserves colimits in the second argument.
\end{opr}

\begin{rem}
\label{DeltaCatEnr}
It is immediate that a $\Delta$-category $\cM$ has a $\SSet$-enrichment given by the mapping spaces
$$
\Map_\cM(x,y)_n := \cM(\Delta^n \otimes x , y).
$$
\end{rem}

\begin{exm}
The terminal category $[0]$ can be equipped with a (trivial) $\Delta$-structure.
\end{exm}

\begin{exm}
\label{DeltaDGMod}
The category $\DGMod_k$ is a $\Delta$-category for $\Delta^n \otimes M := C_\bullet (\Delta^n) \otimes_k M$, where $C_\bullet$ is the chain complex functor. The natural transformation $diag$ comes from the Alexander-Whitney map as follows: 
$$
diag: C_\bullet (\Delta^n) \stackrel{C_\bullet(\delta)}{\to} C_\bullet (\Delta^n \times \Delta^n) \to C_\bullet (\Delta^n) \otimes_k C_\bullet (\Delta^n).
$$ 
\end{exm}

\begin{exm}
Any simplicial model category $\cM$ is a $\Delta$-category in the obvious way.
\end{exm}

\begin{prop}
\label{DeltaLaxAction}
If $\cM$ is a $\Delta$-category then $\otimes : \Delta \times \cM \to \cM$ can be extended uniquely to a functor $\otimes: \SSet \times \cM \to \cM$ such that
\begin{enumerate}
\item $\otimes$ preserves colimits in each argument,
\item there is a family of maps 
\begin{equation}
\label{actionmap}
a(S,T,x):(S \times T) \otimes x \to S \otimes (T \otimes x)
\end{equation}
natural in $S,T \in \SSet$ and $x \in \cM$, associative in a suitable sense and so that for each $[n]$, the composition
$$
\Delta^n \otimes x \to (\Delta^n \times \Delta^n) \otimes x \to \Delta^n \otimes (\Delta^n \otimes x)
$$
equals $diag(n,x)$ of Definition \ref{deltacat}. Moreover, $a(S,T,x)$ is an isomorphism whenever $S$ or $T$ is discrete.
\end{enumerate}
\end{prop}

\noindent We sometimes call the natural map $a(S,T,x)$ the \emph{action map}.

\proof
Recall that to each simplicial set $S$ we can associate its category of simplexes $\Delta/S$. Its objects are all simplexes of $S$, represented as maps $\Delta^n \to S$, and a morphism between two such objects is given by a map $[n] \to [m]$ in $\Delta$ compatible with morphisms to $S$. Let $s: \Delta/S \to \Delta$ denote the functor $(\Delta^n \to S) \mapsto [n]$, and define $S \otimes x := \colim_{\Delta/S} s \otimes x$. 

For two $S, T \in \SSet$, we have a canonical map $\Delta/(S\times T) \to \Delta/S \times \Delta/T$ induced by the two projections $d_s: \Delta/(S\times T) \to \Delta / S$ and $d_t: \Delta/(S\times T) \to \Delta/T$. Denote again by $s: \Delta/S \to \Delta$, $t: \Delta/T \to \Delta$ and also by $st: \Delta/(S \times T) \to \Delta$ the corresponding forgetful functors. Then we have a sequence of maps
$$
(S \times T) \otimes x \cong \colim_{\Delta/(S \times T)} st \otimes x \to \colim_{\Delta/S \times T} st \otimes (st \otimes x) \cong \colim_{\Delta/(S \times T)}  (s \circ d_s) \otimes ((t \circ d_t) \otimes x) \to 
$$
$$
\to \colim_{\Delta_S} \colim_{\Delta_T} s \otimes (t \otimes x) \cong S \otimes (T \otimes x).
$$
Given the construction of this morphism, one can witness the naturality and check its associativity; due to the third condition of Definition \ref{deltacat}, the action map becomes an isomorphism, 
$$
(S \times T) \otimes x \stackrel{\sim}{\to} S \otimes (T \otimes x),
$$
when any of the two $S,T \in \SSet$ is discrete. This proves the last assertion.
\endproof

\begin{exm}
\label{canonical}
For any cocomplete category $\cM$, there is canonical $\Delta$-structure on 
$\Delta^\op \cM = Fun(\Delta^\op,\cM)$, which produces a strict associative action of simplicial sets. Given a simplicial set $K$ and a simplicial object $X \in \Delta^\op \cM$, we define
$$
(K \otimes X)_n = K_n \otimes X_n = \coprod_{K_n} X_n.
$$
\end{exm}

\begin{opr}
	\label{deltafunctordefinition}
Given two categories $\cM, \cN$ with $\Delta$-structures, a \emph{$\Delta$-functor }$F: \cM \to \cN$ is a functor between underlying categories together with a family of morphisms
\begin{equation}
\label{deltafunctorstructuremaps}
m_F([n],x): \Delta^n \otimes F(x) \to F(\Delta^n \otimes x)
\end{equation}
natural in both $[n]$ and $x$. It is required to be compatible with the diagonal maps and unit isomorphisms.
\end{opr}

\begin{rem}
Equivalently, a $\Delta$-functor $F : \cM \to \cN$ is a simplicial functor for the simplicial enrichment mentioned in Remark \ref{DeltaCatEnr}. In particular, it is evident that the composition of $\Delta$-functors is naturally a $\Delta$-functor.

As one can see, the notion of a $\Delta$-functor is lax, in the sense that the arrow (\ref{deltafunctorstructuremaps}) is not invertible, and is in the direction as defined. The reader can compare this notion with that of a lax monoidal functor. One could introduce colax and strict versions of functors between $\Delta$-categories, however, Definition \ref{deltafunctordefinition} is directly related to the examples we have in mind.
\end{rem}

\begin{exm}
\label{DeltaDGModFunctor}
Given two $\Delta$-categories $\cM$ and $\cN$, we adopt the convention that their product $\cM \times \cN$ is equipped with the $\Delta$-structure acting on both components: $\Delta^n \otimes (X,Y) := (\Delta^n \otimes X, \Delta^n \otimes Y)$.  Then, the tensor product of chain complexes 
$$\otimes_k: \DGMod_k \times \DGMod_k \to \DGMod_k,$$ or similarly, for any finite\footnote{$S=*$ corresponds to the identity functor, $S=\emptyset$ corresponds to the inclusion of $k$ in $\DGMod_k$.} set $S$, the $S$-fold tensor product
$$
\otimes_k: \DGMod^S_k \to \DGMod_k
$$
can be naturally equipped with the structure of a $\Delta$-functor. 
\end{exm}

\begin{prop}
A $\Delta$-functor on $F: \cM \to \cN$ between $\Delta$-categories determines a family of maps $m_F(S,x): S \otimes F(x) \to F(S \otimes x)$ natural in $S \in \SSet$ and $x \in \cM$, which restricts to $m_F([n],x)$ for $S = \Delta^n$ and respects the action maps of Proposition \ref{DeltaLaxAction}.
 \end{prop}
\proof Define $m_F(S,x)$ as 
$$
S \otimes F(x) \cong \colim_{\Delta/S} s \otimes F(x) \stackrel{m_F}{\longrightarrow} \colim_{\Delta/S} F(s \otimes x) \to F(\colim_{\Delta/S} s \otimes x) \cong F(S \otimes x). 
$$
Then the result follows. \endproof

Recall that for a functor $F: I^\op \times I \to \cM$, its \emph{coend} \cite[IX.6]{ML} is defined as the universal object $\int^I F$ in $\cM$ together with maps $F(i,i) \to \int^I F$ for each $i \in I$, such that for any morphism $i \to i'$, the induced diagram commutes:
\begin{diagram}[small]
F(i',i) & \rTo & F(i,i) \\
\dTo	&	&	\dTo \\
F(i',i') & \rTo & \int^I F. \\
\end{diagram}
Coends exist in $\cM$ when $\cM$ is cocomplete.

Recall that in a model or higher-categorical setting, simplicial objects correspond to a homotopical version of coequalisers. The corresponding quotient objects are usually given by (homotopy) coends. Accordingly, 

\begin{opr}
Let $\cM$ be a $\Delta$-category, and $X: \Delta^\op \to \cM$ a simplicial object in $\cM$. Its \emph{geometric realisation} is defined as 
$$
|X| := \int^{\Delta^\op} \Delta^\bullet \otimes X
$$
Where $\Delta^\bullet$ is the Yoneda functor (\ref{deltayoneda}). Varying $X$, we get a functor $| - |: \Delta^\op \cM \to \cM$.
\end{opr}

For $S \in \SSet$ and $A \in \cM$, it is evident that the realisation of the simplicial object $[n] \mapsto S_n \otimes A$ is canonically isomorphic to $S \otimes A$.

\begin{prop}
\label{deltastructureassociativity}
For a $\Delta$-functor $f: \cM \to \cN$ we have a canonical natural transformation 
$$
s_f: |f ( -) |_{\cN} \to f | - |_{\cM}  
$$
between the corresponding geometric realisations, where $f : \Delta^\op \cM \to \Delta^\op \cN$ is the induced functor. It is compatible with the composition in the following sense:  the pasting of
\begin{diagram}[small]
\Delta^\op \cM & \rTo^f & \Delta^\op \cN & \rTo^g & \Delta^\op \cK \\
\dTo  &	\stackrel{s_f}{\Downarrow}	& \dTo &	\stackrel{s_g}{\Downarrow}	&  \dTo  \\
\cM			& \rTo_f	&	\cN		&	\rTo_g &	\cK	\\
\end{diagram}
with vertical functors given by realisations, is equal to $s_{gf}$.

\end{prop}
\proof A tedious but straightforward check. \endproof

\subsection{Homotopical $\Delta$-categories}

For any bisimplicial object $X \in (\Delta^{\op} \times \Delta^\op) \cM$ denote by $\delta^* X \in \Delta^\op \cM$ the diagonal simplicial object, that is, the pull-back of $X$ along the diagonal map $\delta: \Delta^\op \to \Delta^\op \times \Delta^\op$.

\begin{opr}
\label{homotopicaldeltacat}
A \emph{homotopical $\Delta$-structure} on a category $\cM$ consists of
\begin{itemize}
    \item a homotopical structure given by a subcategory $\cW \subset \cM$,
    \item a $\Delta$-structure with the functor $\otimes:\Delta \times \cM \to \cM$,
\end{itemize}
so that the following conditions are satisfied:
\begin{enumerate}
    \item the subcategory $\cW$ satisfies two-out-of-six (Definition \ref{twooutofsix}),
    \item $\cM$ is a $\Delta$-category and $\cW$ is preserved by small coproducts,
    \item the induced functor $\otimes : \SSet \times \cM \to \cM$ respects weak equivalences in each variable,
    \item the induced action map (\ref{actionmap}) $a(S,T,x):(S \times T) \otimes x \to S \otimes (T \otimes x)$ is a weak equivalence for each $x \in \cM$ and $S,T \in \SSet$,
    \item the geometric realisation functor $|-|:\Delta^{\op} \cM \to \cM$ preserves pointwise weak equivalences and for each bisimplicial object $X \in (\Delta^{\op} \times \Delta^\op) \cM$, the natural composite map
\begin{equation}
\label{realisationassociativity}
\int^{\Delta^\op} \Delta^\bullet \otimes \delta^* X \to \int^{\Delta^\op} \Delta^\bullet \otimes (\Delta^\bullet \otimes \delta^* X) \to  \int^{\Delta^\op \times \Delta^\op} \Delta^\bullet \otimes (\Delta^\bullet  \otimes X) 
\end{equation}
is a weak equivalence.
\end{enumerate}
A category together with a homotopical $\Delta$-structure is called a \emph{homotopical $\Delta$-category}.
\end{opr}

\begin{rem}
The identity (\ref{realisationassociativity}) implies that for any bisimplicial object $B: \Delta^\op \times \Delta^\op \to \cM$, one may take geometric realisations in any order. See Proposition \ref{propertiesofhomdelta} for details.
\end{rem}

\begin{exm}
\label{SMCat}
Some simplicial model categories $\cM$, for instance simplicial presheaves with injective model structure or simplicial vector spaces, produce examples of homotopical $\Delta$-categories. The non-trivial point here is that the realisation functor $\Delta^\op \cM \to \cM$ only preserves weak equivalences between Reedy cofibrant objects \cite[VII.3.6]{G-J}, but for the model categories just mentioned, all objects of $\Delta^\op \cM$ are automatically cofibrant.
\end{exm}

\begin{exm}
\label{DGMod}
The category $\DGMod_k$ is a homotopical $\Delta$-category for the $\Delta$-structure of Example \ref{DeltaDGMod} and $\cW$ being the class of quasi-isomorphisms. In this case, all simplicial objects are Reedy-cofibrant, and the functor of geometric realisation is known to be left Quillen for the Reedy model structure on simplicial objects \cite[Lemma 9.8]{BER-MOER}. 
\end{exm}

We assemble together some of the properties of geometric realisation. Define the category $\Delta_\infty$ as a subcategory of $\Delta$ consisting of all objects and maps $f:[m] \to [n]$ such that $f(m)=n$. One has the adjunction
$$
j: \Delta \rightleftharpoons \Delta_\infty: i
$$
where $j([n]) = [n+1]$ and the adjunction map $id \to i \circ j$ evaluated on $[n]$ is the inclusion $[n] \hookrightarrow [n+1]$ of $[n]$ as first $n+1$ elements of $[n+1]$. Intuitively, $j$ attaches one more, maximal, element to each $[n]$. 

\begin{opr}
\label{augmentedsimplicialobject}
A \emph{split-augmented} simplicial object is a functor $\bar X: \Delta^{\op}_\infty \to \cM$. A simplicial object $X: \Delta^\op \to \cM$ admits a (split) augmentation iff $X \cong j^* \bar X$ for some $\bar X: \Delta_\infty^\op \to \cM$.
\end{opr}

For a \emph{bisimplicial} object $X: \Delta^\op \times \Delta^\op \to \cM$, we denote by $||X|_2 |_1$ its repeated realisation, that is the coend of the functor
$$
([i],[j],[k],[l]) \mapsto \Delta^i \otimes (\Delta^j  \otimes X_{kl}) 
$$
and by $||X|_1|_2|$ its transpose realisation, which is just a repeated realisation of a transposed bisimplicial object $([n],[m]) \mapsto X_{mn}$.

\begin{prop} 
\label{propertiesofhomdelta}
For a homotopical $\Delta$-category $\cM$, the following is true:
\begin{enumerate}
\item For any simplicial object $X$ admitting an augmentation $\bar X$, its realisation is weakly equivalent to $X_{-1}:= \bar X_0$. Precisely, there are weak equivalences
$$
X_{-1} \to | X | \to X_{-1}
$$
with composition identity that come from the extra maps $X_{-1} \to X_n$ and $X_n \to X_{-1}$. 
\item Given a morphism $X  \to Y $ of bisimplicial objects, we have 
$$
\left( ||X|_2 |_1 \to ||Y|_2 |_1 \right) \in  \cW \Leftrightarrow \left( ||X|_1 |_2 \to ||Y|_1 |_2 \right) \in \cW
$$ 
\end{enumerate}
\end{prop}

The first statement has already been encountered in many different contexts (including higher-categorical ones, see \cite[Lemma 6.1.3.16]{LUHTT}), and is known to largely depend on the combinatorics of split-augmented objects. 
To distinguish simplicial and bisimplicial objects, we write $X_\bullet$ ($[n] \mapsto X_n$) and $X_{\bullet \bullet}$ ($([n],[m]) \mapsto X_{nm})$ for a simplicial and a bisimplicial object, correspondingly.  

\proof  The first statement is proven in a few steps. Denote the tensoring of Example \ref{canonical} by $\circledast: \SSet \times \Delta^\op \cM \to \Delta^\op \cM$. We now prove that the structure of a homotopical $\Delta$-category on $\cM$ gives rise to a family of weak equivalences $|K \circledast X_\bullet| \to K \otimes |X_\bullet|$ natural in $K \in \SSet$ and $X_\bullet \in \Delta^\op \cM$. The maps are constructed as the following sequence:
$$
K \otimes |X_\bullet| = \left( \int^{\Delta^\op} \Delta^\bullet \otimes K_\bullet \right) \otimes \left( \int^{\Delta^\op} \Delta^\bullet \otimes X_\bullet \right) \cong \int^{\Delta^\op \times \Delta^\op} \Delta^\bullet \otimes (\Delta^\bullet \otimes \coprod_{K_\bullet} X_\bullet) 
$$
$$
\leftarrow \int^{\Delta^\op} \Delta^\bullet \otimes \delta^*(\coprod_{K_\bullet} X_\bullet) \cong \int^{\Delta^\op} \Delta^\bullet \otimes (K \circledast X_\bullet) = |K \circledast X_\bullet|. 
$$
Here $\Delta^\bullet \otimes K_\bullet$ denotes the bifunctor $([n],[m]) \mapsto \coprod_{K_m} \Delta^n$. The first line isomorphism is then due to the well-known ``Fubini'' property of coends \cite{ML}. The only non-invertible map in the chain, 
$$
\int^{\Delta^\op \times \Delta^\op} \Delta^\bullet \otimes (\Delta^\bullet \otimes \coprod_{K_\bullet} X_\bullet) \leftarrow \int^{\Delta^\op} \Delta^\bullet \otimes \delta^*(\coprod_{K_\bullet} X_\bullet) 
$$
is a weak equivalence by Definition \ref{homotopicaldeltacat}.

By definition, a simplicial homotopy equivalence in $\Delta^\op \cM$ consists of two maps $f: X_\bullet \to Y_\bullet$ and $g: Y_\bullet \to X_\bullet$, and two diagrams
\begin{diagram}[small,nohug]
X_\bullet					&			&			& & &	Y_\bullet					&			&				\\
\dTo						& \rdTo^{gf}	&			& & &	\dTo						& \rdTo^{fg}	&				\\
\Delta^1 \circledast X_\bullet	& \rTo^h		&	X_\bullet, 	& & &	\Delta^1 \circledast Y_\bullet	& \rTo^{h'}		&	Y_\bullet 		\\
\uTo						& \ruTo_{id}	&			& & &	\uTo						& \ruTo_{id}	&				\\
X_\bullet					&			&			& & &	Y_\bullet					&			&				\\
\end{diagram}
where the vertical maps are induced from the two inclusions $[0] \rightrightarrows [1]$ in $\Delta$. The natural weak equivalence $|K \circledast X_\bullet| \to K \otimes |X_\bullet|$ then implies that, after the realisation, the compositions $|g||f|$ and $|f||g|$ are weak equivalences. By two-out-of-six we get that $|g|$ and $|f|$ are weak equivalences as well.

It is known \cite[Lemma 4.5.1]{RIEHL} that $X_\bullet$ admitting an augmentation $\bar X_\bullet$ in the sense of Definition \ref{augmentedsimplicialobject} leads to a diagram in $\Delta^\op \cM$
$$
\bar X_0 \to X_\bullet \to \bar X_0
$$
naturally appearing from the extra morphisms in $\bar X_\bullet$. The composition of these maps is the identity $id_{X_0}$, and both maps can be shown to be simplicial homotopy equivalences in $\Delta^\op \cM$ for the canonical simplicial structure of Example \ref{canonical}; they thus become weak equivalences after applying geometric realisation.

For the last statement of the proposition, observe that both maps in question are weakly equivalent (that is, weakly equivalent as objects of the category $\cM^{[1]}$ of maps in $\cM$) to the map
$$
\int^{\Delta^\op} \Delta^\bullet \otimes \delta^*(X_{\bullet \bullet}) \to \int^{\Delta^\op} \Delta^\bullet \otimes \delta^*(Y_{\bullet \bullet})
$$ which finishes the proof. \endproof

\section{Fibrations, opfibrations, sections}

The main aim of this section is to provide background material for the language of Grothendieck fibrations (contravariant families of categories) and opfibrations (covariant families). Along with the basic notions, we also explain some operations such as taking powers and transpose fibrations. Some parts of this material are largely folklore with no well-known reference available, so we supply proofs. We finish by considering families of homotopical $\Delta$-categories as described by homotopical $\Delta$-opfibrations.
  
\subsection{Basic notions}

\begin{opr}
Let $p: \cE \to \cC$ be a functor.  A morphism $\alpha: x \to y$ in $\cE$ is \emph{$p$-Cartesian} \cite{SGA1}, or simply Cartesian, if, for every morphism $\beta: z \to y$ of $\cE$ such that $p(\beta)=p(\alpha)$, there exists a unique morphism $\gamma: z \to x$ such that $\beta = \alpha \gamma$ and $p(\gamma)=id_{p(z)}$. 

A morphism $\alpha: x \to y$ in $\cE$ is \emph{$p$-opCartesian} if it is Cartesian for $p^{\op}: \cE^{\op} \to \cC^{\op}$. 
\end{opr}

\begin{opr}
\label{fibrationdefinition}
A functor $p: \cE \to \cC$ is called a \emph{Grothendieck fibration} \cite{SGA1,VIST} (or simply a fibration) of categories iff the following two conditions are satisfied:
\begin{itemize}
\item For every morphism $f: a \to b$ of $\cC$ and $y \in \cE$ such that $p(y)=b$ there exists a Cartesian morphism $\alpha:x \to y$ in $\cE$ covering $\alpha$, that is, $p(\alpha) = f$. 

\item The composition of Cartesian morphisms is a Cartesian morphism.
\end{itemize} 

Dually, $p$ is called an \emph{opfibration} of categories iff $p^{\op}: \cE^{\op} \to \cC^{\op}$ is a fibration of categories. 
\end{opr}

An (op)fibration $q:\cO \to \cC$ is \emph{small} if both $\cC$ and $\cO$ are small.

\begin{const}
\label{grothendieckconstruction}
Given a functor $E$ from $\cC$ to categories, we produce an opfibration, which we denote $\int E \to \cC$ and call the \emph{Grothendieck construction} \cite{VIST} of $E$. An object of $\int E$ is a pair $(c,x)$ of $c \in \cC$ and $x \in E(c)$, and a morphism $(c,x) \to (c',x')$ consists of $f: c \to c'$ together with a map $\alpha: E(f)(x) \to x'$ in $E(x')$. 

In this paper, we already used the symbol $\int$ for coends; in fact, $\int E$ can be reproduced as a certain coend. 

Similar considerations apply in the case of fibrations. For a contravariant category-valued functor $F$ defined on $\cC$, the Grothendieck construction is a fibration $\int F \to \cC$ with same pairs $(c,x)$ serving as objects, but with maps given by $f: c \to c'$ and $\beta: x \to F(f) x'$. 
\end{const}

\begin{exm}
\label{fibsset}
The category $\Delta/S$ for a simplicial set $S$ is exactly the domain of the fibration $\int S \to \Delta$, where we view $S$ as a functor $\Delta^\op \to \Set \subset \Cat$ contravariant on $\Delta$.
\end{exm}

Denote by $\Fin_*$ the category with objects finite sets, and morphisms $f: S \to T$ given by partially defined maps. That is, $f$ is a pair $(D, \tilde f)$ of a subset $D \subset S$ and a map of sets $\tilde f: D \to T$. 

\begin{exm}
\label{monfibstrict}
Consider a strict symmetric monoidal category $\cM$, that is, a category together with an associative commutative unital product $\otimes$.  From this data, we can form a functor $M$ to categories defined on $\Fin_*$ by the assignment $S  \mapsto M(S) := \cM^S.$ For each partially defined map $f: S \to T$, there is a functor $f_!: \cM^S \to \cM^T$, which sends an $\{X_s\}_{s \in S} \in \cM^S$ to $\{Y_t\}_{t \in T} \in \cM^T$ with $Y_t = \otimes_{s \in f^{-1} t} \, X_s.$ When the inverse image is empty, the product is equal to the unit object. 

Now we take the Grothendieck construction of $M$, obtaining $\int M \to \Fin_*$. The category $\cM^\otimes := \int M$ is otherwise described as follows. Its objects are $(S, \{X_s\}_{s \in S})$ where $S \in \Fin_*$ and each $X_s$ is an object of $\cM$. A morphism $(S, \{X_s\}_{s \in S}) \to (T, \{Y_t\}_{t \in T})$ is a partially defined map $f:S \to T$, and a morphism $\otimes_{s \in f^{-1}(t)} X_s \to Y_t$ for each $t \in T$.
\end{exm}

A general (op)fibration is not equal to $\int E \to \cC$ for some functor $E$ to categories. Nonetheless, consider an opfibration $p:\cE \to \cC$. For $c \in \cC$, denote by $\cE(c) = p^{-1}(c)$ the fibre of $p$ over $c$. Let $f: c \to c'$ be a morphism in $\cC$ and $x \in \cE(c)$. Then we can choose an opCartesian morphism $\alpha:x \to y$ such that $p(\alpha) = f$. This specifies an object $y \in \cE(c')$. By the universal property of opCartesian maps, the assignment $x \mapsto y$ defines a functor $f_!:\cE(c) \to \cE(c')$, which is called a transition functor along $f$. One can check \cite{SGA1, VIST} that for each composable pair $f,g$, there exists a coherence isomorphism $g_! \circ f_! \cong (g \circ f)_!$ such that for any composable triple of arrows $f, g, h$, any choice of coherence isomorphisms renders commutative the following diagram:
\begin{diagram}[small]
h_! g_! f_! & \rTo^\sim & h_! (gf)_! \\
\dTo^\sim	&					&	\dTo^\sim \\
(hg)_! f_! & \rTo^\sim & (hgf)_! . \\
\end{diagram}
In the literature, such choice of an assignment $f \mapsto f_!$ together with coherence isomorphisms is called a cleavage. 

\begin{exm}
\label{monfib}
Take an arbitrary symmetric monoidal category $\cM$ Define the category $\cM^\otimes$ with the same objects and morphisms as in Example \ref{monfibstrict}, but now with compositions defined with the help of coherence isomorphisms for $\otimes$. The forgetful functor $\cM^\otimes \to \Fin_*$ is then an opfibration of categories. It is possible to characterise exactly the opfibrations arising from symmetric monoidal categories using Segal conditions \cite{LU, TV}.
\end{exm}

\begin{opr}
Let $p:\cE \to \cC$ and $q: \cE' \to \cC$ be two (op)fibrations. A \emph{lax morphism} between $p$ and $q$ is a functor $F: \cE \to \cE'$ such that $q \circ F = p$. Such $F$ is called a \emph{Cartesian morphism} if it takes (op)Cartesian morphisms of $\cE$ to (op)Cartesian morphisms of $\cE'$. A \emph{section} of an (op)fibration $p$ is a lax morphism from the (op)fibration $id_\cC: \cC \to \cC$ to $p$. 

A morphism between two lax or Cartesian morphisms is a natural transformation $\alpha:F \to F'$ such that for each $x$ in the domain $\cE$, $\alpha_x$ projects to $id_{p(x)}$. 
\end{opr}

We denote by $\Cart(\cE,\cE')$ the category of Cartesian morphisms and by $\Sect(\cC,\cE)$ the category of sections. 

\begin{const}
Take an opfibration $p: \cE \to \cC$, and for each $c \in \cC$, denote by $c \backslash \cC$ the category of objects under $c$ \cite{ML}. The forgetful functor $c \backslash \cC \to \cC$ is an opfibration. Then the assignment $c \mapsto \Cart(c \backslash \cC, \cE)$ defines a covariant category-valued functor on $\cC$. When $\cC$ is small, this construction is inverse up to an equivalence \cite{VIST} to (Grothendieck) Construction \ref{grothendieckconstruction}.
\end{const}

This implies that any opfibration (and, similarly, a fibration) $p: \cE \to \cC$ can be, up to an equivalence, replaced by an opfibration $\tilde p: \tilde \cE \to \cC$, for which the assignment $c \mapsto \cE(c)$ can be made into a strict functor by a choice of transition functors along maps in $\cC$. We call the opfibrations (similarly, fibrations) with later property \emph{strictly cleavable}.

\begin{opr}
In $\Fin_*$, a map $\rho:S \to T$ is \emph{inert} \cite{LU} if it is defined on a subset $T'$ of $S$ isomorphic to $T$, and the restriction $\tilde \rho: T' \to T$ is a bijection. 
\end{opr}

\begin{exm}
\label{exmalgebra}
Any algebra object $A$ in a symmetric monoidal category $\cM$ gives a section $A: \Fin_* \to \cM^\otimes$ by the rule $S \mapsto (A,...,A) \in \cM^\otimes (S)$. Conversely, consider a section $B:\Fin_* \to \cM^\otimes$ which sends the inert maps to opCartesian maps. Then the value $B(1) \in \cM$ of $B$ on a one-element set $1$ becomes a commutative monoid. To show this, note that the `inert-to-opCartesian' condition implies that $B(S) \cong (S, {B(1)}_{s \in S})$ by the means of all opCartesian maps out of $B(S)$ lying over all inert maps $S \to 1$. Next, take the map $S \to 1$ defined on each element of $S$, and consider its image, $B(S) \to B(1)$, in $\cM^\otimes$. It factors as 
$$
B(S) \to B(1)^{\otimes S} \to B(1)
$$
and gives multiplication maps $B(1)^{\otimes S} \to B(1)$. Examining the composition of maps in $\cM^\otimes$, it can be checked that everything is determined when $S$ is a two-element set, and associated operation $B(1) \otimes B(1) \to B(1)$ must be associative, with a unit $I_\cM \to B(1)$ given by the value of $B$ on the map $\emptyset \to 1$ in $\Fin_*$, and commutative. with commutativity following from the action of $B$ on the non-trivial automorphism of $S$. 
\end{exm}

\begin{exm} To explain the term `lax morphism', consider a lax symmetric monoidal functor $F: \cM \to \cN$ between two symmetric monoidal categories. It means that there is a natural family of maps $FX \otimes FY \to F(X \otimes Y)$ and a map $I_\cN \to F I_\cM$ with $I_\cN, I_\cM$ unit objects, which satisfy suitable coherence conditions.  For opfibrations of Example \ref{monfib} the assignment $(S, \{X_s\}_{s \in S}) \mapsto (S, \{FX_s\}_{s \in S})$ then induces a functor $F^\otimes : \cM^\otimes \to \cN^\otimes$ over $\Fin_*$ which does not necessarily preserve Cartesian maps. For example, consider the map $(X,Y) \to X \otimes Y$, which is opCartesian in $\cM^\otimes$. Its image in $\cN^\otimes$, $(FX, FY) \to F(X \otimes Y)$, factors through the $\cN^\otimes$-opCartesian map $(FX, FY) \to FX \otimes FY$ exactly by the lax functor structure of $F$.
\end{exm}

\begin{exm}
\label{laxfunctorgrothendieckconstruction}
Let $L: \int \cE \to \int \cE'$ be a lax morphism between two Grothendieck constructions of covariant functors $\cE, \cE': \cC \to \Cat$. For each $c \in \cC$, $L$ specifies a functor $L_c : \cE(c) \to \cE'(c)$. For each morphism $f: c \to c'$, we get a $2$-square  
\begin{diagram}
 \cE(c)  & \rTo^{L_c} & \cE'(c) \\
\dTo<{\cE(f)}			& \stackrel{L_{f}}{\Leftarrow} & \dTo>{\cE'(f)} \\
\cE(c')  & \rTo^{L_{c'}} & \cE'(c'). \\
\end{diagram}
The natural transformation appears because the image under $L$ of an opCartesian map $X \to \cE(f) X$ ($X \in \cE(c)$) may not be opCartesian. Factoring $LX \to L \cE(f) X$,
$$
LX \to \cE'(f) LX \to L \cE(f) X,
$$
gives $\cE'(f) LX \to L \cE(f) X$; for each $X \in \cE(c)$, all such maps assemble into $L_f$. For two composable arrows $f: c \to c'$, $g: c' \to c''$, there is a pasting property relating $L_f, L_g$ and $L_{g f}$ similar to the one of Proposition \ref{deltastructureassociativity}.

For fibrations, there is a difference on the level of 2-diagrams. Consider $\cF, \cF': \cC^\op \to \Cat$ and take a lax morphism $M: \int \cF \to \int \cF'$ of fibrations over $\cC$. For $f: c \to c'$, we obtain a diagram
\begin{diagram}
\cF(c)  & \rTo^{M_c} & \cF'(c) \\
\uTo<{\cF(f)}			& \stackrel{M_{f}}{\Rightarrow} & \uTo>{\cF'(f)} \\
\cF(c')  & \rTo^{M_{c'}} & \cF'(c') \\ 
\end{diagram}
with $M_f$ given by arrows of the form $M \cF(f) X \to \cF'(f) M X$.
\end{exm}

\begin{opr} 
\label{transposefib}
\label{cospan}
Fix an opfibration $p: \cE \to \cC$. Define a category denoted as $\cE^\top$ as follows:
\begin{enumerate}
\item $Ob(\cE^\top) = Ob(\cE)$ 
\item A morphism from $x \to z$ in $\cE^\top$ is an isomorphism class of cospans in $\cE$
$$
x \longrightarrow y \longleftarrow z
$$
such that the left arrow is fibrewise, $p(x \to y) = id_{p(x)}$, and the right arrow is opCartesian. 
\end{enumerate} 
There is an evident functor $p^\top: \cE^\top \to \cC^\op$ which sends maps $x \longrightarrow y \longleftarrow z$ to $p(y \longleftarrow z)$. A morphism of $\cE^\top$ is $p^\top$-Cartesian iff it can be represented by a span of the form $y \stackrel{id_y}{\longrightarrow} y \longleftarrow z$.  The functor $p^\top$ is a fibration, which we call the \emph{transpose fibration} of $p$.
\end{opr}

If $\cE \to \cC$ equals $\int E \to \cC$ for a functor $E: \cC \to \Cat$, then $\cE^\top \to \cC^\op$ is equivalent to the (fibrational) Grothendieck construction applied to $E: (\cC^\op)^\op \to \Cat$ viewed as a contravariant functor on $\cC^\op$.  

\begin{exm}
\label{monoidaldualfib}
The transpose fibration $(\cM^\otimes)^\top \to \Fin_*^\op$ of $\cM^\otimes \to \Fin_*$ from Example \ref{monfib} can be constructed by hand just like the original opfibration. Now, to define a lax functor between two such fibrations, $(\cM^\otimes)^\top \to (\cN^\otimes)^\top$, we would have to consider \emph{oplax} symmetric monoidal functors between symmetric monoidal categories. Such functors preserve (cocommutative) coalgebra objects, as opposed to lax symmetric monoidal functors, which preserve commutative algebra objects. In particular, any coalgebra object $C$ of $\cM$ produces a section of $(\cM^\otimes)^\top \to \Fin_*^\op$ by the assignment $S \mapsto (S, \{ C \}_S)$.    
\end{exm}

Given a functor $F:\cD \to \cC$, we can pull back (op)fibrations over $\cC$ to $\cD$, with the result again being (op)fibrations. Similarly, given a section $A: \cC \to \cE$ of an (op)fibration $\cE \to \cC$ we obtain from it the section $F^* A: \cD \to F^* \cE$ of the pull-back (op)fibration $F^* \cE \to \cD$. This operation defines the pull-back functor $F^*: \Sect(\cC,\cE) \to \Sect(\cD,F^*\cE)$. 

\begin{prop}
\label{2functoriality}
Assume given a \emph{fibration} $\cF \to \cC$ and a natural transformation $\alpha: F \to G$ of functors $F,G:\cD \to \cC$. Then
\begin{itemize} 
\item there is a natural Cartesian map of fibrations $R_\alpha: G^* \cF \to F^* \cF$, which we call the \emph{restriction} map,
\item given a section $A: \cC \to \cF$, there is a natural morphism of sections
$$
F^* A \to R_\alpha G^* A.
$$
\end{itemize}
\end{prop}
The fact that $\cF \to \cC$ is a fibration, and not an opfibration, is important for the direction of the arrows in the proposition.

\proof Up to an equivalence we can assume $\cF \to \cC$ to be strictly cleavable. Take $d \in \cD$. For each object $X$ of $\cE(G(d)) = G^* \cF (d)$, we have a Cartesian arrow $Y \to X$ in $\cF$ over $\alpha_d: F(d) \to G(d)$. The value $R_\alpha X$ is then defined to be equal to $Y$; its action on morphisms can be defined similarly.   

Given a section $A$, its value on $\alpha_d : F(d) \to G(d)$ can be naturally factored as
$$
F^*A(d) = A(F(d)) \to R_\alpha A(G(d)) \to A(G(d)) = G^* A(d).
$$
Varying $d$, the arrows $F^* A(d) \to R_\alpha A(G(d)) = (R\alpha G^*A ) (d)$ define the natural transformation in question. \endproof

\begin{opr}
\label{PowerFib2}
\label{poweringexplicit}
Let $p:\cE \to \cC$ be an (op)fibration and $I \in \Cat$ a category.
\begin{itemize}
\item A product of $I$ and $p: \cE \to \cC$ is the functor $I \times p: I \times \cE \to \cC, \, \, \, \, \, (i,x) \mapsto p(x).$ 
\item A powering of $p$ with $I$ is the functor $p^I: \cE^I \to \cC$ where $\cE^I$ is the subcategory of $Fun(I,\cE)$ consisting of all functors $F: I \to \cE$ such that $p \circ F$ is a constant functor $I \to \cC$. 
\end{itemize}
Both these functors are (op)fibrations. 
\end{opr}

Unfortunately, the choice of notation such as $\cE^I$ may lead to confusion as before we used it to denote the whole category of functors $I \to \cE$. We thus adopt a convention that for (op)fibrations, the powering notation works in the sense of definition above and not otherwise\footnote{If we think of ordinary categories as Grothendieck fibrations over a point, then there is no notational ambiguity.}.

\begin{opr}
\label{powerfib}
\label{powerfibconstr}
Given a fibration $p:\cF \to \cC$ and an opfibration $q:\cO \to \cC$ with small fibres, a \emph{power fibration} $p^q: \cF^\cO \to \cC$ is defined as follows. An object of $\cF^\cO$ is a pair of $c \in \cC$ and a functor $X: \cO(c) \to \cF$ such that $p X$ is constant of value $c$. A morphism $(c,X) \to (c',Y)$ consists of $f: c \to c'$ and a natural transformation $X \to Y \circ f_!$ of functors $\cO(c) \to \cF$ for some choice of transition functor $f_!: \cO(c) \to \cO(c')$. The functor $\cF^\cO \to \cC$ is the natural projection.
\end{opr}

One can verify that $\cF^\cO \to \cC$ is again a fibration, with fibres equivalent to $Fun(\cO(c),\cF(c))$. A fibrational transition functor $\cF^\cO(c') \to \cF^\cO(c)$ is given by precomposing an object $F: \cO(c) \to \cF(c)$ with $f_!: \cO(c) \to \cO(c')$ and postcomposing with $f^*: \cF(c') \to \cF(c)$ for some choice of transition functors $f_!$ and $f^*$ in $\cO$ and $\cF$ respectively.   

\begin{lemma}
\label{powerfibprop}
For a functor $F: \cD \to \cC$, and $p:\cF \to \cC, q: \cO \to \cC$ as above, 
\begin{enumerate}
\item There is an equivalence of categories
$$
\Sect(\cO, q^* \cF) \cong \Sect (\cC, \cF^\cO).
$$
\item There is a Cartesian map
$$
(F^* \cF)^{F^* \cO} \to F^* (\cF^\cO)
$$
which is moreover an equivalence over $\cC$.  
\end{enumerate}
\end{lemma}
\proof Clear. \endproof

\begin{prop}
\label{Fibpushforward}
Let $p: \cF \to \cC$ be a fibration with cocomplete fibres.
\begin{enumerate}
\item For any functor $X: I^\op \times I \to \cF$ such that $p X$ is constant of value $c \in \cC$, its coend $\int^I X$ exists in $\cF$ and can be calculated in $\cF(c)$, defining a lax morphism $\int^I : \cF^{I^\op \times I} \to \cF$ of fibrations over $\cC$.
\item Let
\begin{diagram}[small]
\cO &		        & \rTo^P	&		& I \times \cC \\
 		     &	\rdTo	&     		& \ldTo		& \\
		     &		        & \cC	&		&  \\
\end{diagram}
be an opCartesian morphism of small strictly cleavable opfibrations. Then the obvious functor
$$
P^*: \Sect(\cC, \cF^I) \to \Sect(\cC,\cF^\cO)
$$
admits a left adjoint $P_!$.
\end{enumerate}
\end{prop}
The fact that $\cF \to \cC$ is a fibration is important in this proposition.

\proof Let $X : I^\op \times I \to \cF$ be such that $p X$ is constant at $c \in \cC$. Take its coend in $\cF(c)$, and denote it $\int^I_c X$. We need to check that it satisfies the universal property of a coend in $\cF$. Let $Z \in \cF$ be an object together with maps $X(i,i) \to Z$ such that for any morphism $i \to i'$, the diagram below commutes:
\begin{diagram}[small]
X(i',i) & \rTo & X(i,i) \\
\dTo	&	&	\dTo \\
X(i',i') & \rTo & Z. \\
\end{diagram}
Applying $p: \cF \to \cC$ to diagrams like above, we find that for each $i \in I$, the map $X(i, i) \to Z$ lies over a fixed morphism of $\cC$ with domain $c$; we denote it as $f: c \to c'$. The fact that $p$ is a fibration then implies that there exists an object $f^*Z$ over $c$ and a Cartesian map $\alpha:f^*Z \to Z$ covering $f$, so that each map $X(i,i) \to Z$ factors through $\alpha$ and the resulting squares
\begin{diagram}[small]
X(i',i) & \rTo & X(i,i) \\
\dTo	&	&	\dTo \\
X(i',i') & \rTo & f^*Z \\
\end{diagram}
are all commutative. This implies the existence of a map $\int_c^I X \to f^*Z$. One can then see that the composition $\int_c^I X \to f^*Z \to Z$ is unique and independent of the choice of $f^* Z \to Z$ made above. One can see then that this map is the one which ensures the universality of $\int_c^I X$ in the whole of $\cF$.
 
We now turn to the second statement. We define a functor on the level of fibrations, $P_!: \cF^\cO \to \cF^I$, by setting $P_! X$ to be the left Kan extension \cite{ML} of $X: \cO(c) \to \cF(c)$ along the map $P_c: \cO(c) \to I$. Following the line of argumentation as for the coends before, we can see that this Kan extension exists and moreover can be calculated in the fibre $\cF(c)$. Given a section $S: \cC \to \cF^\cO$, we apply $P_!$ to its values, inducing the sought-after functor $P_!$, left adjoint to the pull-back $P^*:\Sect(\cC, \cF^I) \to \Sect(\cC,\cF^\cO)$. \endproof

\subsection{Homotopical $\Delta$-opfibrations}

\begin{opr}
\label{homotopicalopfibration}
A \emph{homotopical structure} on an opfibration $\cE \to \cC$ consists of a homotopical structure on $\cE$, given by a subcategory $\cW \subset \cE$ of weak equivalences, compatible with the opfibration in the following sense: 
\begin{enumerate}
\item the image of $\cW$ in $\cC$ consists of identity morphisms,
\item in a commutative square
\begin{diagram}[small]
A & \rTo^\alpha & B \\
\dTo^{f} & & 	\dTo>{f'} \\
A' & \rTo^{\alpha'} & B' \\
\end{diagram}
if we have $f \in \cW$ and $\alpha, \alpha'$ are opCartesian, then $f' \in  \cW$.
\end{enumerate}
\end{opr}

The definition of a homotopical fibration is dual, with implication in the diagram above going into opposite direction. One can see that for an (op)fibration with homotopical structure, $\cW = \coprod_{c \in \cC} \cW(c)$, with $(\cE(c),\cW(c))$ being a homotopical category for each $c \in \cC$. The transition functors of the (op)fibration send $\cW(c)$ to $\cW(c')$. 

\begin{opr}
	\label{wcart}
Given an opfibration $\cE \to \cC$ with a homotopical structure, a morphism $\alpha: x \to y$ of $\cE$ is \emph{weakly opCartesian} if it can be factored as an opCartesian morphism $x \to \alpha_! x$ followed by a weak equivalence $\alpha_! x \to y$.
Dually, one has the notion of a weakly Cartesian morphism for homotopical fibrations, where the order of factorisation is reversed.
\end{opr}

Before proceeding with the $\Delta$-structure, we need some extra notation. Take an opfibration $\cE \to \cC$ and consider its transpose fibration (Definition \ref{transposefib}) $\cE^\top \to \cC^\op$. Now, take the product fibration $\Delta \times \cE^\top \to \cC^\op$. If we have any functor $\otimes: \Delta \times \cE^\top \to \cE^\top$ over $\cC^\op$, then, for a map $f: c \to c'$ of the original base category $\cC$, we get, after relabelling $\cE^\top(c) = \cE(c)$ and $\cE^\top(f) = f_!: \cE(c) \to \cE(c')$, the following diagram
\begin{diagram}
\Delta \times\cE(c)  & \rTo^{\otimes_c} & \cE(c) \\
\dTo<{id \times f_!}			& \stackrel{m_{f}}{\Rightarrow} & \dTo>{f_!} \\
\Delta \times \cE(c')  & \rTo^{\otimes_{c'}} & \cE(c') \\
\end{diagram}
so that each $m_f: - \otimes_{c'} f_! - \to f_!(- \otimes_c -)$  appears in the same way as in Example \ref{laxfunctorgrothendieckconstruction}. Moreover, $f \mapsto m_f$ is suitably functorial in $f$.

\begin{opr}
\label{deltastructurefib}
A $\Delta$-structure on an \emph{fibration} $\cF \to \cC$ is a $\Delta$-structure $\otimes: \Delta \times \cF \to \cF$ such that
\begin{enumerate}
\item $\otimes$ is a lax morphism of fibrations over $\cC$,

\item the natural transformation $diag$ and unitality isomorphism (see Definition \ref{deltacat}) are fibrewise.

\item For each $f: c \to c'$  of $\cC$ and a choice of $f^*: \cF(c') \to \cF(c)$ via the fibration property, the induced diagram 
\begin{diagram}
\Delta \times\cF(c)  & \rTo^{\otimes_c} & \cF(c) \\
\uTo<{id \times f^*}			& \stackrel{m_{f}}{\Rightarrow} & \uTo>{f^*} \\
\Delta \times \cF(c')  & \rTo^{\otimes_{c'}} & \cF(c') \\
\end{diagram}
makes $f^*$ into a $\Delta$-functor.
\end{enumerate}

A $\Delta$-structure on an \emph{opfibration} $\cE \to \cC$ consists of a $\Delta$-structure on its transpose fibration. In particular, the notion is not dual.
\end{opr} 
 
\begin{exm} The prototype example for how one ought to think about such $\Delta$-structures is the following.  Consider a covariant functor $F$ on $\cC$ such that each $F(c)$ is a $\Delta$-category and each $F(c \to c')$ is a $\Delta$-functor. Forgetting the $\Delta$-structure, we can apply the Grothendieck construction, obtaining an opfibration $\int F \to \cC$. This opfibration is, then, equipped with a $\Delta$-structure inherited from the values of $F$ on the objects and morphisms of $\cC$.  
\end{exm}

\begin{opr}
	\label{homdeltaopfdef}
A \emph{homotopical} $\Delta$-(op)fibration is an (op)fibration $\cE \to \cC$ together with a homotopical structure and a $\Delta$-structure, such that for each $c \in \cC$, the induced structure on the fibre $\cE(c)$ is that of a homotopical $\Delta$-category.
\end{opr}

\begin{exm}
Chain complexes give us homotopical $\Delta$-opfibration $\DGMod^\otimes \to \Fin_*$. The $\Delta$-structure on the opfibration is essentially explained in Examples \ref{DeltaDGMod} and \ref{DeltaDGModFunctor}, and the weak equivalences are simply induced from the quasi-isomorphisms of $\DGMod_k$. In the same way, those simplicial model categories which give us a homotopical $\Delta$-structure (Example \ref{SMCat}) can as well give us homotopical $\Delta$-opfibrations. If such a category $\cM$ in addition possesses a compatible monoidal structure (this is true, for example, both for simplicial presheaves and for simplicial vector spaces), then the associated opfibration $\cM^\otimes \to \Fin_*$ is a homotopical $\Delta$-opfibration.  
\end{exm}

\begin{prop}
\label{fibrationrealisation}
Let $\cF \to \cC$ be a homotopical $\Delta$-fibration. Then
\begin{itemize}
\item For any functor $F: \cD \to \cC$, the pull-back $F^* \cF \to \cD$ is again a homotopical $\Delta$-fibration.
\item There is a lax \emph{realisation} morphism of fibrations 
\begin{diagram}[small]
\cF^{\Delta^\op}	&		& \rTo^{|-|} 	&			&	\cF \\
	& \rdTo	&				&	\ldTo	&		\\
	&		&	\cC		&			&	\\
\end{diagram}
such that on each fibre, the functor $\Delta^\op \cF(c) \to \cF(c)$ is the geometric realisation for the $\Delta$-category $\cF(c)$.
\end{itemize}
\end{prop}
\proof Recall that $\cF^{\Delta^\op}$ is the subcategory of $Fun(\Delta^\op,\cF)$ consisting of $X : \Delta^\op \to \cF$ which become constant after composing with $\cF \to \cC$. Since $\cF$ has a $\Delta$-structure, we can consider the assignment $X \mapsto \Delta^\bullet \otimes X$, which defines a lax morphism $\cF^{\Delta^\op} \to \cF^{\Delta \times \Delta^\op}$ of fibrations over $\cC$ (it can be seen to not preserve Cartesian arrows). We then use the first part of Proposition \ref{Fibpushforward} and compose the just-obtained functor with the coend $\cF^{\Delta \times \Delta^\op} \to \cF$ to obtain the realisation functor of this proposition. \endproof

\begin{rem}
\label{algcoalg}
From the perspective of $\Delta$-structures, we see that it is preferable to consider fibrations and not opfibrations. The motivation for opfibrations as basic ingredients of the play comes from Example \ref{exmalgebra}, where algebra objects in a monoidal category $\cM$ are presented as sections of $\cM^\otimes \to \Fin_*$. Suitably normalised sections of the transpose fibration $(\cM^\otimes)^\top \to \Fin_*^\op$ correspond, on the other hand, to coalgebra objects in $\cM$, as can be seen from Example \ref{monoidaldualfib}. The fact that formalism of the subsequent chapter presents derived sections of an opfibration as certain sections of the related transpose fibration may be perceived as an example of certain ``Bar-Cobar'' duality relating algebraic and coalgebraic objects.
\end{rem}

\section{Derived sections}

This section is devoted to our definition of derived sections. Starting from the introduction of simplicial replacements of categories, we assemble together a few homotopy-theoretic facts to be used later. We also outline how one prolongs an opfibration $\cE \to \cC$ to the simplicial replacement $\bC$ of $\cC$. We finish by defining presections and derived sections. 

As shown in \cite{BK} and later in \cite{RIEHL}, simplicial replacements of categories allow calculating homotopy colimits and may be viewed as certain cofibrant replacements in homotopy-algebraic sense. We will make this perspective on simplicial replacements concrete in our setting through the constructions of the next section. 

\subsection{Simplicial replacements}

\begin{opr} \label{simplicialreplacement} Given a small category $\cC$, its \emph{simplicial replacement}, denoted $\bC$, is the \emph{opposite} of $\Delta/N \cC = \int N \cC$, that is the opposite of the category of simplexes of the simplicial set $N \cC$ (cf. Example \ref{fibsset}).
\end{opr}

An object of $\bC$ is a sequence $c_0 \to ... \to c_n$ of composable morphisms in $\cC$. Any functor $F: \cD \to \cC$ induces a functor $\bF:\bD \to \bC:$ by the rule $\bF(d_0 \to ... \to d_n)= Fd_0 \to ... \to Fd_n$. Observe that $\bF$ commutes with the projections from $\bD$ and $\bC$ to $\Delta^\op$. 

Our choice of terminology is inspired by the simplicial replacement of diagrams in \cite{BK}. Note that $\bC$ is an ordinary category and not a simplicial category. The assignment $\cC \mapsto \bC$ defines a functor from $\Cat$ to the full subcategory of $\Cat/(\Delta^\op)$, consisting of opfibrations over $\Delta^\op$ with discrete fibres. 

\begin{ntn}
\label{joinobject}
We often denote by $\pi: \bC \to \Delta^\op$ the natural projection. An object $c_0 \to ... \to c_n$  of $\bC$ will be denoted as $\bc_{[n]}$ (so that $\pi(\bc_{[n]})=[n]$) or simply as $\bc$ when its underlying $\Delta$-object is not important. Given two objects $\bc_{[n]}$, $\bc'_{[m]}$, and a map $\alpha: c_n \to c'_0$, we denote by $\bc_{[n]} *^\alpha \bc'_{[m]}$ the 'concatenated' object
$$
c_0 \to ... \to c_n \stackrel{\alpha}{\to} c_0' \to ... \to c_n'.
$$
\end{ntn}

\begin{lemma}
\label{headandtail}
There are functors $h_\cC:\bC \to \cC$ and $t_\cC: \bC \to \cC^\op$ given by $\bc_{[n]} \mapsto c_0$ or $\bc_{[n]} \mapsto  c_n$ respectively. \endproof
\end{lemma}

\begin{opr}
A map $\zeta:\bc_{[n]} \to \bc'_{[m]}$ is \emph{anchor} iff its projection in $\Delta$, $\pi(\zeta):[m] \to [n]$, is an interval inclusion of $[m]$ as first $m+1$ elements of $[n]$, i.e. $\pi(\zeta)(i) = i$ for $0 \leq i \leq m$. In particular, $m$ should be less or equal than $n$.

A map $\zeta:\bc_{[n]} \to \bc'_{[m]}$ is \emph{structural} iff its image under $t_\cC$ is an identity and the underlying map in $\Delta^\op$ preserves the endpoints: $\pi(\zeta)(m) = n$.

We denote by $A_\bC$ and $S_\bC$ the sets of all anchor and structure maps respectively.
\end{opr}

The prototypical example to think of is the span diagram \begin{diagram}[small,nohug]
         &  		&	c_0 \to ... \to c_n &				&	\\
		& \ldTo		& 							&	\rdTo	&		\\
c_0 \to ... \to c_k &	 		 & 								&				&	c_n	\\
\end{diagram}
with $k \leq n$, the left map being anchor and the right map being structural. In general, structural maps can correspond to non-injections in $\Delta$.

\begin{prop}[Factorisation and localisation properties] 
\label{locprop}
For a small category $\cC$,
\begin{enumerate}
\item Every map $\bc \to \bc'$ can be uniquely factored as an anchor map $\bc \to \bc''$ followed by a structural map $\bc'' \to \bc$. 

\item Any functor $F: \bC \to \cN$ which sends anchor maps of $\bC$ to isomorphisms factors essentially uniquely as $F = \tilde F \circ h_\cC$ for $\tilde F: \cC \to \cN$. In other words, $\cC$ is a localisation of $\bC$ with respect to anchor maps.
\end{enumerate}

\end{prop}

\proof The factorisation property is clear and is inherited from $\Delta$. For the second statement, we first note that the functor $h_\cC^*: \cN^\cC \to \cN^\bC$ is full and faithful (see e.g. \cite[Section 4.4]{RIEHL}). Let $F: \bC \to \cN$ be a functor which sends $A_\bC$ to isomorphisms. Define a new functor $\bar F: \cC \to \cN$. On objects, $\bar F(c) = F(c)$ where $c$ is viewed as an object of $\bC$ of zero length. Take a span 
$$
c  \longleftarrow (c \stackrel f \to c' ) \longrightarrow c',
$$
the action of $F$ on it gives a span $F(c) \leftarrow F(c \to c') \to F(c')$. Inverting the left arrow we get a map $\bar F(f): \bar F(c) \to \bar F(c')$. The action of $F$ on objects of higher length, $c \to c' \to c''$, and on degenerate objects, $c \stackrel{id}{\to} c$, then ensures that $\bar F$ is indeed a functor and $F = \bar F h_\cC$.\endproof

\begin{rem}
To stress, the class of anchor maps is not saturated in the sense one applies when one speaks of localisation \cite{DHKS}. Indeed, not every map which becomes an isomorphism under $h_\cC$ is an anchor map. 
\end{rem}

The proposition permits to justify the idea that, given a homotopical category $\cM$ with weak equivalences $\cW$, a functor $F: \bC \to \cM$ sending $A_\bC$ to $\cW$ is a suitable weakening of the concept of a functor from $\cC$ to $\cM$. The action of $F$ on spans in $\bC$ like 
$$
c  \longleftarrow (c \stackrel f \to c' ) \longrightarrow c',
$$
where the left arrow is an anchor map, gives a span $F(c) \stackrel \cW \leftarrow F(c \to c') \rightarrow F(c')$, where the left map is a weak equivalence. On the level of $\Ho \cM$, this span gives a map $F(c) \to F(c')$, which one can denote $F(f)$. Applying $F$ to higher-length objects then ensures higher coherences for the `weak functor' $F$.

The spans of the form $X \stackrel \cW \leftarrow Y \rightarrow Z$ have appeared before many times in the context of localisation (for example, they are known under the name 'cocycles' in \cite{JAR}). For an arbitrary homotopical category $\cM$, such spans may not constitute a good presentation of morphisms in $\Ho \cM$. In practice one may need to make additional assumptions about $\cM$. For our purposes, homotopical $\Delta$-categories provide a sufficient answer. 

\bigskip

We conclude this subsection by summing up a few results concerning the functors from a simplicial replacement. Let $I$ be a small category and denote by $\mathbb I$ its simplicial replacement. 

\begin{opr}
\label{realisationofafunctor}
Let $\cM$ be a homotopical $\Delta$-category. For $X: \mathbb I \to \cM$, its \emph{realisation} is defined as $| \Pi X|$, where $|-|$ is the geometric realisation for $\cM$ and $\Pi: Fun(\bI, \cM) \to \Delta^\op \cM$ is the left Kan extension along the canonical projection $\pi: \mathbb I \to \Delta^\op$.
\end{opr}

For any object $i \in I$, there is naturally a map $X(i) \to |\Pi X|$. The following lemmas illustrate that the behaviour of the assignment $X \mapsto |\Pi X|$ are resemblant of that of a homotopy colimit of $X: \bI \to \cM$, thought of as a ``weak functor'' from $I$ to $\cM$. 

\begin{lemma} 
\label{terminalobject}
Let $I$ be a category with a terminal object $1$, and $\cM$ be a homotopical $\Delta$-category. Then any $X: \mathbb I \to \cM$ sending the anchor maps $A_{\mathbb I}$ to maps in $\cW$, the natural map $X(1) \to | \Pi X|$ is an equivalence.
\end{lemma}
\proof Consider an 'augmented' functor $X^{aug}: \mathbf{i} \mapsto X(\mathbf i  *^{x} 1)$ (here $x$ corresponds to the canonical map to the terminal object $t_I(\mathbf{i}) \to 1$). It is then easy to see that there is a canonical equivalence $X^{aug} \to X$ coming from the maps $X(\mathbf i *^{x} 1) \to X(\mathbf i)$. It then becomes an equivalence of realisations. The object $\Pi X^{aug}$, however, can be completed to a split-augmented simplicial object $\tilde X^{aug}:\Delta_\infty^\op \to \cM$ defined by the formula
$$
\tilde X^{aug}_n = \Pi X^{aug}_{n-1}, \, \, \, n>0,
$$
$$
\tilde X^{aug}_0 = X(1).
$$
in particular, one augmentation map $X(1) \to \tilde X^{aug}_1 = \coprod_{i} X(i \to 1)$ comes from the image $X(1) \to X(1 \to 1)$ of the degeneracy map $1 \to ( 1 \to 1)$ and the other map
$$
\tilde X^{aug}_1 = \coprod_{i} X(i \to 1) \to X(1)
$$  
is just the coproduct of the natural maps $X(i \to 1) \to X(1)$. By Proposition \ref{propertiesofhomdelta} we have the equivalences
$$
X(1) \to |\Pi  X^{aug}| \to X(1)
$$
and we can see that the composite map $X(1) \to |\Pi X^{aug}| \to |\Pi X|$, which is an equivalence, is equal to the map in question. \endproof

\begin{lemma}
\label{hcdiag}
Let $I$ be a category with contractible nerve and $\cM$ be a homotopical $\Delta$-category. If a functor $X: \mathbb I \to \cM$ takes all morphisms of $\mathbb I$ to isomorphisms, then the natural map $X(i) \to | \Pi X |$ is an equivalence for any $i \in I$.
\end{lemma}
\proof Fix $i \in I$. Proposition \ref{locprop} implies that the functor $X$ can be factored as $\overline X \circ h_{I}$ with $\overline X: I \to \cM$. $X$ moreover factors through the fundamental groupoid of $I$, which is contractible. One can then see that
$$
\Pi X_n = \coprod_{\mathbf{i}_{[n]}} X(\mathbf{i}_{[n]}) \cong \coprod_{\mathbf{i}_{[n]}} (\overline X \circ h_{I})(\mathbf{i}_{[n]}) \cong \coprod_{\mathbf{i}_{[n]}} X(i_0) \cong \coprod_{\mathbf{i}_{[n]}} X(i),
$$
and so $|\Pi X| = NI \otimes X(i)$, which is equivalent to $X(i)$, and the map $X(i) \to |\Pi X|$ in question is a homotopy inverse of the projection $NI \otimes X(i) \to X(i)$. \endproof

\subsection{Homotopical category of derived sections}

We now study how opfibrations over $\cC$ interact with the simplicial replacement $\bC$.

\begin{opr}
	\label{opfibsimpextdef}
For an opfibration $\cE \to \cC$, its \emph{simplicial extension} is a \emph{fibration} $\bfE \to \bC$ which is the pull-back of a transpose fibration $\cE^\top \to \cC^\op$ along $t_\cC: \bC \to \cC^{\op}$.
\end{opr}

\begin{rem}
We stress that $\bfE$ is not a simplicial replacement of $\cE$ or $\cE^\top$. In particular, the fibre of $\bfE \to \bC$ over an object $\bc_{[n]}$ is equivalent to $\cE(c_n)$. If $\cE \to \cC$ comes from a functor $\cE: \cC \to \Cat$, then $\bfE \to \bC$ corresponds to the functor 
\begin{diagram}
\bC^\op & \rTo^{t^\op_\cC} & \cC & \rTo^{\cE} & \Cat
\end{diagram}
viewed as a contravariant functor on $\bC$. 
\end{rem}

\begin{rem}
\label{fibrationnormalisation}
Given two functors $k_1, k_2: K \to \bC$ and a natural transformation $\alpha: k_1 \rightarrow k_2$ valued in structural maps $S_\bC$, we have that the induced Cartesian map of fibrations 
$$
\alpha^*: k^*_2 \bfE \to k^*_1 \bfE
$$
is in fact an equivalence.
\end{rem}

We can also pull back $\cE \to \cC$ to $\bC$ by the means of the functor $h_\cC: \bC \to \cC$. The link between this pull-back and the fibration $\bfE \to \bC$ is in the following: 

\begin{prop}
\label{EtoEtop}
Given an opfibration $p: \cE \to \cC$, there is a morphism $T:h^*_\cC \cE \to \bfE$ commuting with functors to $\bC$ which sends opCartesian maps of $h^*_\cC \cE$ to Cartesian maps of $\bfE$ and is universal, i.e. any other functor $G: h^*_\cC \cE \to \bfE$ over $\bC$ with such a property factors through $T$ up to a natural isomorphism.
\end{prop}
\proof Consider the category $\cX$ defined as follows.
\begin{itemize}
\item An object of $\cX$ is a pair $(\bc_{[n]}, \alpha)$ where $\bc_{[n]} = c_0 \to ... \to c_n$ is an object of $\bC$ and $\alpha: x \to y$ is an opCartesian map in $\cE$ which covers the composition $c_0 \to c_n$ in $\cC$ (i.e. $p(\alpha) = c_0 \to c_n$),
\item A morphism $(\bc_{[n]}, \alpha:x \to y) \to (\bc'_{[m]}, \beta: x' \to y')$ consists of a map $\bc \to \bc'$ in $\bC$ and a map $\gamma: x \to x'$ which covers the induced map $c_0 \to c'_0$.
\end{itemize} 
One can check that the natural functor $\cX \to \bC$ is an opfibration, and it is easy to see that the assignment $(\bc, \alpha: x \to y) \mapsto (\bc, x)$ defines an equivalence over $\bC$ of opfibrations $\cX \stackrel \sim \to h_\cC^* \cE$.

On the other hand, consider the assignment $(\bc, \alpha: x \to y) \mapsto (\bc, y)$. We claim that it defines a functor $\bar T: \cX \to \bfE$ commuting with projections to $\bC$. Let $(f,t):(\bc, \alpha: x \to y) \to (\bc', \beta: x' \to y')$ be a map. In particular, we have the following diagram in $\cE$:
\begin{equation}
\label{diagramfort}
\begin{diagram}[small,nohug]
x 						& \rTo^t 	& x'	 				\\
\dTo<\alpha		&				&	\dTo>\beta \\
y						&				& 	y'.		\\
\end{diagram}
\end{equation}
Suppose first that the map $t$ is fibrewise. Then by opCartesian property there exists a map $t': y \to y'$ rendering the diagram (\ref{diagramfort}) commutative. Remembering the description of arrows in Definition \ref{cospan}, we define $\bar T(f,t) = (f,y \stackrel{t'}{\to} y' \stackrel{id}{\leftarrow} y')$; in other words, we view $t'$ as a fibrewise map of $\cE^\top$ 

Next, if $t$ is opCartesian, find an opCartesian map $k: y' \to z$ in $\cE$ covering $c'_m \to c_n$ (which is induced from $f:\bc \to \bc'$). The composition $k \beta t$ and $\alpha$ both project along $\cE \to \cC$ to the map $c_0 \to c_n = c_0 \to c'_0 \to c'_m \to c_n$, hence there is a (fibrewise) isomorphism $z \cong y$. This implies that the diagram (\ref{diagramfort}) can be completed as  
\begin{diagram}[small,nohug]
x 						& \rTo^t 	& x'	 				\\
\dTo<\alpha		&				&	\dTo>\beta \\
y						&	\lTo^{t'}			& 	y'.	\\
\end{diagram}
with all arrows opCartesian in $\cE$. We put, again, $\bar T(f,t) = (f,y \stackrel{id}{\to} y \stackrel{t'}{\leftarrow} y')$, thus viewing $t'$ as a Cartesian map of $\cE^\top$. Any other case of $(f,t)$ can be treated by reducing to these two cases. 

Inverting the equivalence $\cX \stackrel \sim \to h_\cC^* \cE$ and composing with $\bar T$, we obtain the desired functor $T: h_\cC^* \cE \to \bfE$, and one can verify its universal property.
\endproof

\begin{opr}
	\label{psectdef}
Given an opfibration $\cE \to \cC$, its category of \emph{presections} is the category
$$
\PSect(\cC,\cE):= \Sect_\bC(\bC, \bfE).
$$
\end{opr}

Recall the functors $h_\cC$ and $T$ discussed before in Lemma \ref{headandtail} and Proposition \ref{EtoEtop}.

\begin{prop}
\label{sectinpsect}
The assignment $S \mapsto T \circ (h_\cC^* S)$ defines a functor $i:\Sect(\cC,\cE) \to \PSect(\cC,\cE)$. Its essential image consists of the presections sending the anchor maps $A_\bC$ to Cartesian morphisms in $\bfE$.
\end{prop}
\proof Note that for any anchor map $a:\bc_{[n]} \to \bc_{[k]}$ a map in $h_\cC^* \cE$ is opCartesian over $a$ iff it is an isomorphism $x \stackrel \sim \to x$ in $\cE(c_0)$. On one hand, the functor $T$ sends such maps to Cartesian maps in $\bfE$; on the other hand, the pull-back section $h^*_\cC S: \bC \to h^*_\cC \cE$ sends anchor maps $A_\bC$ precisely to identities in $\cE$. Further details are then clear. \endproof

\begin{rem}
Consider an object $\bc_{[n]} = c_0 \stackrel{f_1}{\longrightarrow} c_1 \stackrel{f_2}{\longrightarrow} ... \stackrel{f_n}{\longrightarrow} c_n$ of $\bC$. Then $S \in Sect(\cC, \cE)$ is sent by the functor above to $i(S)$ such that $i(S)(\bc_{[n]}) \cong (f_n ... f_1)_! S(c_0)$ where $(f_n ... f_1)_!: \cE(c_0) \to \cE(c_n) = \bfE(\bc_{[n]})$ is a transition functor along the composition of $f_i$.
\end{rem}

Assume now that $\cE \to \cC$ has a homotopical structure $\cW$.

\begin{opr}
The \emph{standard homotopical structure} on $\PSect(\cC,\cE)$ is defined by the subcategory of those morphisms $A \to A'$ for which the map $A(\bc_{[n]}) \to A'(\bc_{[n]})$ is in $\cW$ for each $\bc_{[n]} \in \bC$.
\end{opr}

We henceforth assume this homotopical structure whenever dealing with $\PSect(\cC,\cE)$. We denote by $\Ho \PSect(\cC,\cE)$ the corresponding localisation.

\begin{opr}
	\label{dsectdef}
A presection $A: \bC \to \bfE$ is a \emph{derived section} iff $A$ sends anchor maps to weakly Cartesian (Definition \ref{wcart}) morphisms in $\bfE$.
\end{opr}

We denote by $\DSect(\cC,\cE)$ the full subcategory of $\PSect(\cC,\cE)$ spanned by derived sections. We restrict the standard homotopical structure from $\PSect(\cC,\cE)$ to $\DSect(\cC,\cE)$ and denote by $\Ho \DSect(\cC,\cE)$ the corresponding localisation. The notion of a derived section depends only on a homotopical structure of $\cE \to \cC$, and not a choice of a $\Delta$-structure.

\section{The pushforward functor}

Given a functor $F: \cD \to \cC$, there is an induced pull-back morphism
$$
\bF^*: \PSect(\cC,\cE) = \Sect(\bC,\bfE) \to \Sect(\bD,F^* \bfE) = \PSect (\cD,F^*\cE)
$$
which restricts well to
$$
\bF^*: \DSect(\cC,\cE) \to \DSect(\cD,F^* \cE)
$$
and is moreover homotopical. It is natural to ask if such a functor may admit a homotopy adjoint. 

In this section, we provide a partial answer to this question, by constructing 
$$
\bF_!: \PSect(\cD,F^* \cE) \to \PSect(\cC, \cE),
$$ 
a homotopical functor in the other direction, together with natural spans relating $\bF_! \bF^*$ and $\bF^* \bF_!$ with identity functors (Proposition \ref{propositionadjunction}). The functor $\bF_!$ may be viewed as an almost left adjoint to the pull-back functor $\bF^*$. If $\bF_!$ preserves derived sections, the spans mentioned above indeed give well-defined maps $\bF_! \bF^* \to id$ and $id \to \bF^* \bF_!$ on the level of homotopy categories of derived sections, satisfying a triangle identity. Even under a weaker assumption that the span relating $\bF_! \bF^*$ to the identity functor consists of weak equivalences, Corollary \ref{ffcondition} implies that $\bF^*$ is full and faithful on homotopy level.

\subsection{Main construction}

Recall that, given two functors $\cD \stackrel F \to \cC \stackrel G \leftarrow \cB$, one can form the corresponding \emph{comma category} $F / G$ \cite{ML}. Its objects are triples $(d,b,\alpha: F(d) \to G(b))$ for $d \in \cD$ and $b \in \cB$. We shall need the following adaptation of this notion:

\begin{opr}
	\label{FovrG}
	Given two a diagram $\cD \stackrel F \to \cC \stackrel G \leftarrow \cB$, the associated \emph{simplicial comma object} $\bF \ovr \bG$ is defined as the \emph{opposite} of the category $\int F \ovr G$, where $F \ovr G: \Delta^\op \times \Delta^\op \to \Set$ is the bisimplicial set
	$$
	F \ovr G([n],[m]) = \{ \bd_{[n]}, \mathbf b_{[m]}, \alpha: F(d_n) \to G(b_0) \}
	$$
	viewed as a contravariant functor to $\Cat$. In other words, it is the category with objects given by triples $(\bd_{[n]}, \mathbf b_{[m]}, \alpha: F(d_n) \to G(b_0))$ and a map
	$$
	(\bd_{[n]}, \mathbf b_{[m]}, \alpha: F(d_n) \to G(b_0)) \to (\bd'_{[k]}, \mathbf b'_{[l]}, \beta: F(d'_k) \to G(b'_0))
	$$
	consists of two maps $\bd_{[n]} \to \bd'_{[k]}$, $\mathbf b_{[m]} \to \mathbf b'_{[l]}$ such that the induced square commutes:
	\begin{diagram}[small,nohug]
	F(d_n) & \rTo^\alpha & G(b_0) \\
	\uTo	&			 & \dTo		\\
	F(d'_k) &	\rTo^\beta & G(b'_0). \\
	\end{diagram}
\end{opr}
We often write $\bD \ovr \bG$ or $\bF \ovr \bC$ instead of $\bF \ovr \bG$ if $F$ or $G$ is the identity functor. Given an object $c \in \cC$, we also consider $\bF \ovr c$ where we treat $c$ as a functor $[0] \to \cC$ and denote its simplicial replacement by the same letter. The canonical functor $\bF \ovr \bG \to \Delta^\op \times \Delta^\op$ is an opfibration with discrete fibres $\bF \ovr \bG([n],[m]) = F \ovr G ([n],[m])$. 

There is a \emph{concatenation} functor $con:\Delta \times \Delta \to \Delta$, $([n],[m]) \mapsto [n] * [m] = [n+m+1]$, and we think that $[n]$ is included as first $n+1$ elements of $[n+m+1]$ and $[m]$ as last $[m+1]$ elements. The action of $con$ on morphisms is then evident.

Then we observe the following. There is a diagram in $\Cat$
\begin{equation}
\label{FovrGdiagram}
\begin{diagram}[nohug]
&       & \bF \ovr \bG &       &     \\
& \ldTo^{pr_\bD} & \Leftarrow \dTo \Rightarrow	& \rdTo^{pr_\bB} &      \\
\bD & \rTo^{\bF}  & \bC                   & \lTo^{\bG}  &  \bB \\
\end{diagram}
\end{equation}
with the middle map denoted $pr_{\bF \ovr \bG}$, covering the diagram
\begin{equation}
\label{DeltaDeltadiagram}
\begin{diagram}[nohug]
&       & \Delta^\op \times \Delta^\op &       &     \\
& \ldTo^{\pi_1} & \Leftarrow \dTo \Rightarrow	& \rdTo^{\pi_2} &      \\
\Delta^\op & \rTo^{id}  & \Delta^\op                   & \lTo^{id}  &  \Delta^\op \\
\end{diagram}
\end{equation}
with the middle map acting as $([n],[m]) \mapsto [n]*[m]$. Moreover,
\begin{itemize}
	\item the left natural transformation $pr_{\bF \ovr \bG} \to \bF \circ pr_\bD$ is valued in anchor maps $A_\bC$,
	\item the right natural transformation $pr_{\bF \ovr \bG} \to \bG \circ pr_\bB$ is valued in structural maps $S_\bC$, 
	\item $pr_{\bB}$ is an opfibration whose classifying functor $\bB \to \Cat$ sends anchor maps to equivalences of categories.
\end{itemize}

All this is evident from Definition \ref{FovrG}: $pr_\bD$ maps $(\bd_{[n]}, \mathbf b_{[m]}, \alpha: F(d_n) \to G(b_0) )$ to $\bd_{[n]}$, $pr_\bB$ maps it to $\mathbf b_{[m]}$, and $pr_{\bF \ovr \bG}$ maps it to $\bF(\bd_{[n]})*^\alpha \bG(\mathbf b_{[m]})$. 

The fact that $pr_{\bB}$ sends anchor maps to equivalences and Proposition \ref{locprop} suggests that that $pr_{\bB}$ can be obtained as a pull-back of an opfibration over $\cB$ along the first element map $h_\cB: \bB \to \cB$. This opfibration $X \to \cB$ consists of the category $X$ whose objects are triples $(\bd_{[n]}, b, \alpha: F(d_n) \to G(b) )$ and morphisms are given by $\bd_{[n]} \to \bd'_{[m]}$ in $\bD$, $b \to b'$ in $\cB$, such that the square 
\begin{diagram}[small,nohug]
F(d_n) & \rTo^\alpha & G(b) \\
\uTo	&			 & \dTo		\\
F(d'_m) &	\rTo^\beta & G(b') 
\end{diagram}
commutes in $\cB$, and the functor $X \to \cB$ is the projection $(\bd_{[n]}, b, \alpha) \mapsto b$.

\begin{lemma}
	\label{homopfibtofib}
	Let $\cE \to \cC$ be a homotopical $\Delta$-opfibration. Then the simplicial extension $\bfE \to \bC$ (Definition \ref{opfibsimpextdef}) is a homotopical $\Delta$-fibration.
\end{lemma}
\proof Clear, as $\bfE \to \bC$ is a pull-back of $\cE^\top \to \cC$, which is a homotopical $\Delta$-fibration by definition (see Definitions \ref{deltastructurefib} and \ref{homdeltaopfdef}). \endproof

We are now ready to outline the construction. Fix a functor $F: \cD \to \cC$ and a homotopical $\Delta$-opfibration $\cE \to \cC$. Consider the diagram (\ref{FovrGdiagram}) for $G = id_\cC$:
\begin{equation}
\label{FovrCdiag}
\begin{diagram}[nohug]
      &       & \bF \ovr \bC &       &     \\
      & \ldTo^{pr_\bD} & \Leftarrow \dTo \Rightarrow	& \rdTo^{pr_\bC} &      \\
\bD & \rTo^{\bF}  & \bC                   & \lTo^{id_\bC}  &  \bC. \\
\end{diagram}
\end{equation}
The middle map is $pr_{\bF \ovr \bC}$. This diagram gives us in particular the restriction morphism of Proposition \ref{2functoriality}
\begin{equation}
\label{restrictionmorph}
R_\bF: (\bF pr_\bD)^* \bfE \to pr^*_{\bF \ovr \bC} \bfE.
\end{equation}
This is a map of fibrations over $\bF \ovr \bC$. 

Next, we observe that there are equivalences
\begin{equation}
\label{descenteqspan}
\Sect(\bF \ovr \bC, pr^*_{\bF \ovr \bC} \bfE) \stackrel{\sim}{\leftarrow} \Sect(\bF \ovr \bC, pr^*_{\bC} \bfE) \stackrel{\sim}{\to} \Sect( \bC,  \bfE^{\bF \ovr \bC}) 
\end{equation}
where the right equivalence is provided by the first assessment of Lemma \ref{powerfibprop} (keep in mind that $pr_\bC$ is a small opfibration). The left map comes from the equivalence
$$
pr^*_{\bC} \bfE \stackrel{\sim}{\to} pr^*_{\bF \ovr \bC} \bfE
$$
provided by Remark \ref{fibrationnormalisation}. We denote by 
\begin{equation}
\label{descenteq}
D_\bF:\Sect(\bF \ovr \bC, pr^*_{\bF \ovr \bC} \bfE) \stackrel{\sim}{\to} \Sect( \bC,  \bfE^{\bF \ovr \bC}) 
\end{equation}
the resulting equivalence constructed from (\ref{descenteqspan}).

There is a natural 'projection' functor $\Pi_\bF$ over $\bC$,
\begin{diagram}[small,nohug]
\bF \ovr \bC &		        & \rTo^{\Pi_\bF}	&		& \Delta^\op \times \bC \\
 		     &	\rdTo<{pr_\bC}	&     		& \ldTo		& \\
		     &		        & \bC	&		&  \\
\end{diagram}
which acts as $(\bd_{[n]}, \bc_{[m]},\alpha: F(d_{n}) \to c_0) \mapsto ([n],\bc_{[m]})$. Exponentiating and taking sections, we obtain a functor $\Pi_\bF^*: \Sect(\bC,\bfE^{\Delta^\op}) \to \Sect(\bC, \bfE^{\bF \ovr \bC})$.

\begin{prop}
\label{pidirect}
The functor
$$
\Pi_\bF^*: \Sect(\bC,\bfE^{\Delta^\op}) \to \Sect(\bC, \bfE^{\bF \ovr \bC})
$$
admits a homotopical left adjoint
\begin{equation}
\label{fromovrtosimp}
\Pi_{\bF,!}:\Sect(\bC, \bfE^{\bF \ovr \bC}) \to \Sect(\bC,\bfE^{\Delta^\op})
\end{equation}
\end{prop}
\proof
See Proposition \ref{Fibpushforward} for the construction of $\Pi_{\bF,!}$. To observe that it is homotopical, note that for each $\bc$, the functor $\bF \ovr \bC(\bc) \to \Delta^\op$ is a \emph{discrete} opfibration, and the pushforward along it amounts to taking coproducts, which are homotopical. 
\endproof

Take a $\cD$-presection $S:\bD \to \bF^* \bfE$. Then apply functors (\ref{restrictionmorph}), (\ref{descenteq}) and (\ref{fromovrtosimp}) to obtain
\begin{equation}
\label{sectionbeforepi}
B_\bullet (S) := \Pi_{\bF,!} D_\bF ( R_\bF \circ pr^*_\bD S) \in Sect(\bC, \bfE^{\Delta}).
\end{equation}

Lemma \ref{homopfibtofib} implies that $\bfE \to \bC$ is a homotopical $\Delta$-fibration. Applying the realisation functor $|-|$ from Proposition \ref{fibrationrealisation}, we get the following:

\begin{opr}
The \emph{derived pushforward} of a presection $A: \bD \to \bF^* \bfE$ is defined as
$$
\bF_! (S) := |B_\bullet (S) | = | \Pi_{\bF,!} D_\bF ( R_\bF \circ pr^*_\bD S) |.
$$
this defines a homotopical functor $\bF_!: \PSect(\cD,\cE) \to \PSect(\cC,\cE)$.
\end{opr}

Since restriction functor and the equivalence $D_\bF$ preserve weak equivalences, in conjunction with Definition \ref{homotopicaldeltacat} and Proposition \ref{pidirect} we indeed have that the functor $\bF_!$ is homotopical.

\begin{rem}
Over an object $\bc_{[m]} = c_0 \stackrel{f_1}{\to} ... \stackrel{f_m}{\to} c_m$, we have 
$$
B_{n}(S)(\bc_{[m]}) = \coprod_{\bd_{[n]}, \alpha: F(d_n) \to c_0}(f_m ... f_1 \alpha)_! S(\bd_{[n]})
$$
where $(f_m ... f_1 \alpha)_!$ is the transition functor $\cE(F(d_n)) \to \cE(c_m)$. This expression is very similar to the bar construction (cf. \cite{BK,RIEHL}). The reason for the fact that coproducts and not more complex colimits appear in the expression is because the fibres of $\bF \ovr \bC(\bc_{[m]}) \to \Delta^\op$ are sets, and that is very similar to the classical case of \cite{RIEHL}. The value $\bF_! S(\bc_{[m]})$ is then the realisation of the simplicial object $B_{n}(S)(\bc_{[m]})$. 
\end{rem}

\subsection{Unit and counit correspondences}

Given a $\cC$-presection $A: \bC \to \bfE$, use $pr_{\bF \ovr \bC}$ from the diagram (\ref{FovrCdiag}) and functors (\ref{restrictionmorph}), (\ref{descenteq}) and (\ref{fromovrtosimp}) to obtain

\begin{equation}
\label{sectionaugbeforepi}
B^{\bF}_\bullet(A) := \Pi_{\bF,!} D_\bF ( pr^*_{\bF \ovr \bC}A) \in Sect(\bC, \bfE^{\Delta}).
\end{equation}
Denote by $A^{\bF}$ the realisation of $B_\bullet^{\bF}(A)$.

\begin{rem}
Again, one can see that explicitly 
$$
B^\bF_{n}(A)(\bc_{[m]}) = \coprod_{\bd_{[n]}, \alpha: F(d_n) \to c_0} A(\bF(\bd_{[n]}) *^\alpha \bc_{[m]}).
$$

\end{rem}

\begin{lemma}
\label{counit}
There is a natural (in $A$) correspondence in $\PSect(\cC,\cE)$
$$
\bF_! \bF^* A \leftarrow A^{\bF} \rightarrow A
$$
coming from the realisation of the correspondence of simplicial presections
$$
B_\bullet(\bF^* A) \leftarrow B_\bullet^{\bF}(A) \rightarrow A
$$
where the rightmost term is a constant simplicial object. When $A$ is a derived section, the left morphisms in the correspondences above are weak equivalences.
\end{lemma}
\proof
First, the construction. Given a $\cC$-presection $A: \bC \to \bfE$, Proposition \ref{2functoriality} and the left triangle of the diagram (\ref{FovrCdiag}) gives us a map of $pr_{\bF \ovr \bC}^* \bF^* \bfE$-sections over $\bF \ovr \bC$
$$
pr^*_{\bF \ovr \bC} A \to R_\bF (\bF pr_\bD)^* A
$$
which is an equivalence when $A$ is a derived section. Indeed, over an object $(\bd,\bc,\alpha)$ of $\bF \ovr \bC$ the map looks like
$$
A(\bF(\bd) *^\alpha \bc) \to (f_n...f_1 \alpha)_! A(\bF(\bd))
$$
with $\bc = c_0 \stackrel{f_1}{\to} ... \stackrel{f_n}{\to} c_n$, and this map is an equivalence precisely because of the derived section condition for $A$. Applying the equivalence $D_\bF$ of (\ref{descenteq}) and then $\Pi_!$ of (\ref{fromovrtosimp}), we get the map
$$
B_\bullet^{\bF}(A) \to B_\bullet(\bF^* A)
$$
between (\ref{sectionbeforepi}) and (\ref{sectionaugbeforepi}) which is again a weak equivalence when $A$ is a derived section.

Proposition \ref{2functoriality} and the right triangle of the diagram (\ref{FovrCdiag}) give us a map
$$
pr^*_{\bF\ovr \bC} A \to pr^*_{\bC} A
$$
and we again apply $\Pi_! D_\bF$. Observe that $\Pi_! D_\bF pr^*_\bC A$ is the following simplicial presection: 
$$
(\Pi_! D_\bF pr^*_\bC A)_n(\bc) = \coprod_{\bd_{[n]}, \alpha: F(d_n) \to c_0} A(\bc) \cong N (F/c_0)(n) \otimes A(\bc).
$$
There is thus a natural map $\Pi_! D_\bF pr^*_\bC A \to A$ to the constant simplicial presection $A$.

The realisation of $\Pi_! D_\bF pr^*_\bC A$ is the presection given by the assignment $\bc \mapsto N (F/c_0) \otimes A(\bc)$. On this level as well, we get the map 
$$
A^\bF \to N (F/h_\cC(-)) \otimes A \to A
$$
which completes the construction. \endproof

\begin{lemma}
\label{degeneraciestoequivalences}
Let $A: \bC \to \bfE$ be a derived section and $s:\bc_{[n]} \to \bc'_{[m]}$ be such a map in $\bC$ that its underlying map $s:[m] \to [n]$ in $\Delta$ is a surjection.
Then $A(s)$ is weakly Cartesian in $\bfE$. 
\end{lemma}
\proof Both $A(\bc_{[n]})$ and $A(\bc'_{[m]})$ are seen to land in the same category, and are weakly equivalent to $f_! A(c_0)$, where $f_! :c_0 \to c_n$ is the map given by taking the composition of arrows in $\bc_{[n]}$. These weak equivalences are moreover compatible with $A(s)$. \endproof 

\begin{lemma}
\label{counitidentity}
For $F = id_\cC$ and a derived section $A: \bC \to \bfE$ both morphisms in the span 
$$
id_{\bC !} id_\bC^* A \leftarrow A^{id_\bC} \rightarrow A
$$
of Lemma \ref{counit} are weak equivalences.
\end{lemma}
\proof Fix $\bc \in \bC$. In the case of the identity functor, we see that $B_\bullet^{id_\bC} A(\bc)$ can calculated as the realisation (cf. Definition \ref{realisationofafunctor}) of the functor $X:\bC \ovr c_0 \to \bfE(\bc)$ defined by the assignment
$$
X((\bc'_{[k]}, \alpha:c'_k  \to c_0)) = A( \bc'_{[k]}*^\alpha \bc).
$$
The category $\bC \ovr c_0$ is the simplicial replacement of the category $\cC/c_0$, and the latter has a terminal object. By Lemma \ref{terminalobject}, the natural map $A(c_0 *^{id_{c_0}} \bc) = X(c_0) \to |\Pi_! X| = A^{id_\bC}(\bc)$ is an equivalence. 

There is also an equivalence $A(\bc) \to A(c_0 *^{id_{c_0}} \bc)$ which comes from the degeneracy $\bc \to c_0 *^{id_{c_0}} \bc$ (cf. Lemma \ref{degeneraciestoequivalences}). One can then see that the composition
$$
A(\bc) \to A(c_0 *^{id_{c_0}} \bc) \to A^{id_\bC} (\bc) \to A(\bc)
$$
is the identity (it is such already on the level of corresponding simplicial objects; also note that the composition $\bc \to c_0 *^{id_{c_0}} \bc \to \bc$ in $\bC$ is the identity $id_\bc$). Thus the $\bc$-th component of the map $A^{id_\bC} \to A$ is an equivalence as a right inverse of an equivalence $A(\bc) \to A(c_0 *^{id_{c_0}} \bc) \to A^{id_\bC} (\bc)$. \endproof

\begin{lemma}
\label{corrlemma3}
For a functor $F: \cD \to \cC$ and a $\bD$-presection $A$, there is a natural (in $A$) morphism
\begin{diagram}[small]
id_{\bD !} id_\bD^* A & \rTo & \bF^* \bF_! A.
\end{diagram}
\end{lemma}
\proof By definition, $id_{\bD !} id_\bD^* A$ is the realisation of the simplicially valued presection $X$ which at $d_0 \stackrel{g_1}{\to} ... \stackrel{g_m}{\to} d_m$ takes the value
$$
[n] \mapsto X_n= \underset{\alpha: d'_n \to d_0}{\coprod_{d'_0 \to ... \to d'_n}}(F(g_m ... g_1 \alpha))_! A(\bd'_{[n]}).
$$
In the case when we calculate $\bF^* \bF_! A$ at $d_0 \stackrel{g_1}{\to} ... \stackrel{g_m}{\to} d_m$, we have the following simplicial object $Y$:
$$
[n] \mapsto Y_n = \underset{\beta: F(d'_n) \to F(d_0)}{\coprod_{d'_0 \to ... \to d'_n}}(F(g_m ... g_1) \beta)_! A(\bd'_{[n]}).
$$
The assignment of $\alpha: d'_n \to d_0$ to $F \alpha: F(d'_n) \to F(d_0) $ induces the map of sets
\begin{equation}
\label{indexingsets}
\{ d'_0 \to ... \to d'_n, \alpha: d'_n \to d_0 \} \to \{ d'_0 \to ... \to d'_n, \beta: F(d'_n) \to F(d_0) \}
\end{equation}
and we obtain a map $X_n \to Y_n$ as $X_n$ and $Y_n$ are the coproducts indexed by the sets in (\ref{indexingsets}). Varying $[n] \in \Delta$, we assemble a map $X \to Y$ of simplicial objects, which after realisations gives the map in question, $id_{\bD !} id_\bD^* A \longrightarrow \bF^* \bF_! A$. \endproof

We finally prove the main proposition of this section:
\begin{prop}
\label{propositionadjunction}
Let $F : \cD \to \cC$, $A \in \DSect(\cC,\cE)$ and $R \in \DSect(\cD, F^* \cE)$.
\begin{enumerate}
\item There is a natural (in $A$) span of presections
\begin{equation}
\label{corrcounit}
\bF_! \bF^* A \leftarrow A^{\bF} \rightarrow A
\end{equation}
which induces a natural transformation $\epsilon: \bF_! \bF^*  \to id$ of functors $\Ho \DSect(\cC,\cE) \to \Ho \PSect(\cC,\cE)$ (where $id$ is the inclusion functor).

\item There is a natural (in $R$) sequence of morphisms
\begin{equation}
\label{corrunit}
R \leftarrow R^{id_\bD} \rightarrow id_{\bD !} id^*_{\bD} R \rightarrow \bF^* \bF_! R
\end{equation}
which induces a natural transformation $\eta : id \to \bF^* \bF_!$ of functors $\Ho \DSect(\cD, F^* \cE) \to \Ho \PSect(\cD, F^* \cE)$.
\item \emph{(Triangle identity)} For each $A \in \Ho \DSect(\cC,\cE)$, the composition in $\Ho \PSect(\cD, F^* \cE)$
\begin{equation}
\label{corrtriangle}
\begin{diagram}
\bF^* A & \rTo^{\eta_{\bF^* A} } & \bF^* \bF_! \bF^* A & \rTo^{\bF^* \epsilon_A} & \bF^* A
\end{diagram}
\end{equation}
is the identity.
\end{enumerate}
\end{prop}
\proof
We proved the first two claims in the preceding lemmas. Only the triangle identity remains. Using the correspondences obtained before, we write a string of morphisms
$$
\bF^*A \stackrel{\sim}{\leftarrow} (\bF^* A)^{id_\bD} \stackrel{\sim}{\rightarrow} id_{\bD !} id^*_{\bD} \bF^* A \rightarrow \bF^* \bF_! \bF^* A \stackrel{\sim}{\leftarrow} \bF^* (A^\bF) \rightarrow \bF^* A 
$$
with all the weak equivalences drawn as $\stackrel{\sim}{\rightarrow}$ or $\stackrel{\sim}{\leftarrow}$.  We can redraw this sequence, obtaining the (potentially non-commutative) diagram
\begin{diagram}[small]
		&			& (\bF^* A)^{id_\bD}	& \rTo^\sim	& id_{\bD !} id^*_{\bD} \bF^* A \\
		& \ldTo^\sim	&					&			&						\\
\bF^*A	&			&					&			&	\dTo					\\
		& \luTo		&					&			&						\\
		&			& \bF^* (A^\bF	)		& \rTo^\sim	&	\bF^* \bF_! \bF^* A		\\
\end{diagram}
The third claim is then equivalent to the commutativity of this diagram. We proceed as follows:  writing down in components the simplicial object used to obtain $(\bF^* A)^{id_\bD}$, we see
$$
(\bF^* A)^{id_\bD} \longleftrightarrow B^{id_\bD}_{n}(\bF^* A)(\bd_{[m]}) = \coprod_{\bd'_{[n]}, \alpha: d'_n \to d_0} A(\bF(\bd'_{[n]} *^\alpha \bd_{[m]})).
$$
In the same way,
$$
\bF^*(A^{\bF}) \longleftrightarrow (\bF^* B^{\bF}_{n}( A))(\bd_{[m]}) = \coprod_{\bd'_{[n]}, \beta: F(d'_n) \to F(d_0)} A(\bF(\bd'_{[n]}) *^\beta \bF(\bd_{[m]})).
$$
Assigning $\alpha \mapsto F(\alpha)$, we see that there is a natural in $A$ map $(\bF^* A)^{id_\bD} \to \bF^* (A^\bF)$. Moreover, a comparison with the construction of Lemma \ref{corrlemma3} reveals that in the resulting diagram
\begin{diagram}[small]
		&			& (\bF^* A)^{id_\bD}	& \rTo^\sim	& id_{\bD !} id^*_{\bD} \bF^* A \\
		& \ldTo^\sim	&					&			&						\\
\bF^*A	&			&	\dTo				&			&	\dTo					\\
		& \luTo		&					&			&						\\
		&			& \bF^* (A^\bF	)		& \rTo^\sim	&	\bF^* \bF_! \bF^* A		\\
\end{diagram}
both the left-hand triangle and the right-hand square commute. \endproof

\begin{corr}
\label{ffcondition}
Assume that for a functor $F: \cD \to \cC$, both maps in the span (\ref{corrcounit}) are weak equivalences. Then $\bF^*: \Ho \DSect(\cC,\cE) \to \Ho \DSect (\cD,F^*\cE)$ is full and faithful.
\end{corr}

\proof This result can be proven as a particular case of the following categorical result:

Let $f:\cM \rightleftarrows \cN:u$ be two functors, $i:\cN_0 \subset \cN$ and $j: \cM_0 \subset \cM$ are full subcategories such that $u   i$ is contained in $\cM_0$. In other words, there is a functor $u_0: \cN_0 \to \cM_0$ with $u i = j u_0$. Suppose furthermore that there are natural transformations $\epsilon: f u i \stackrel \sim \to i$ and $\eta: j \to uf j$ defined over $\cN_0$ and $\cM_0$ respectively such that the triangle identity is satisfied: $ui = ju_0 \to ufj u_0 = ufui \to ui$ is the identity. Then $u_0$, or equivalently $u i$, is full and faithful.  

In turn, the categorical result is proven as follows. The functoriality of $u_0$ supplies us with maps $u(x,y):\cN_0(x,y) \to \cM_0(u_0 x, u_0 y)=\cM(uix, uiy)$. Given a map $\alpha: u_0 x \to u_0 y$, we define $v(x,y) \alpha$ to be the map fitting in the commutative square
\begin{diagram}[small]
fu i x & \rTo^{f \alpha} & fu i y \\
\dTo<{\epsilon_x}>\sim &			& \dTo>{\epsilon_y}<\sim \\
i x &	\rTo_{v(x,y) \alpha} &	i y \\
\end{diagram}
(here we use that $i$ is a full and faithful inclusion). This defines the map $v(x,y): \cM(uix, uiy) \to \cN_0(x,y)$ which is inverse to $u(x,y)$.
\endproof

Note in particular that in the situation like above, for $A \in \DSect(\cC,\cE)$, $\bF_! \bF^* A$ is again a derived section.

\section{Case of a resolution}

In this section, we study the functors of $\bF^*$ and $\bF_!$ for a particular class of functors $F$, which we call resolutions:

\begin{opr}
\label{catresolution}
A functor $F: \cD \to \cC$ is a \emph{resolution} if it is an opfibration and each fibre $\cD(c)$ is contractible (that is, its nerve $N \cD(c)$ is contractible). 
\end{opr}


Denote by $\bD(c)$ the simplicial replacement of $\cD(c)$.

\begin{opr}
\label{locconst}
Let $F: \cD \to \cC$ be a resolution. A presection $A: \bD \to \bF^* \bfE$ is \emph{locally constant} if for any fibre $\cD(c)$ over $c \in \cC$, the composite functor
$$
\bD(c) \to \bD \stackrel{A}{\to} \bfE(c) = \cE(c) 
$$
sends all morphisms of the domain to weak equivalences. A derived presection is locally constant if it is locally constant as a presection.
\end{opr}

We denote by $\PSect(\cD,\cE)_{lc}$ and $\DSect(\cD,\cE)_{lc}$ the corresponding full homotopical subcategories of locally constant (pre)sections. It is clear that any (pre)section of the form $\bF^* A$ is locally constant.

When $F$ is a resolution, we can prove two general results concerning $\bF^*$. Here is the first result:

\begin{thm}
\label{maintheorem1}
Let $\cE \to \cC$ be a homotopical $\Delta$-opfibration and $F: \cD \to \cC$ be a resolution (Definition \ref{catresolution}). Then after passing to localisations, the pull-back functor $\bF^*:\Ho \DSect(\cC , \cE) \to \Ho \DSect(\cD, \cE)$ is full and faithful.
\end{thm}
The proof of this theorem relies on the machinery of pushforwards considered in the previous chapter leading to Corollary \ref{ffcondition} and the additional manipulations are similar in the spirit to the proof of Cofinality Theorem in \cite{BK}. Indeed, we shall prove that $\bF_! S(\bc_{[m]})$ which is calculated as the realisation of 
$$
[n] \mapsto B_{n}(S)(\bc_{[m]}) = \coprod_{\bd_{[n]} \subset \bD, \alpha: F(d_n) \to c_0}(f_m ... f_1 \alpha)_! S(\bd_{[n]})
$$
can be also calculated, up to a coherent zigzag of equivalences, as the realisation of  
$$
[n] \mapsto  \coprod_{\bd_{[n]} \subset \bD(c_0) }(f_m ... f_1 )_! S(\bd_{[n]}).
$$
When $S = \bF^* A$ for a derived section $A: \bC \to \bfE$, the second realisation is seen to be equivalent to $A(\bc_{[m]})$. The difficulty of the proof comes mostly from the complexity of objects involved, and the necessity to make sure that the aforementioned zigzag of equivalences is (equivalent, for a fixed $\bc$, to) the counit correspondence (\ref{corrcounit}) of Proposition \ref{propositionadjunction}.   

We can also characterise the homotopical essential image of $\bF^*$. Unfortunately, we only know how to do it for $F$-special (cf. Definition \ref{special}) homotopical $\Delta$-opfibrations: 
\begin{thm}
\label{maintheorem2}
Let $F: \cD \to \cC$ be a resolution and $\cE \to \cC$ be a $F$-special homotopical $\Delta$-opfibration. Then the functor $\bF^*:\Ho \DSect(\cC , \cE) \to \Ho \DSect(\cD, \cE)_{lc}$ is an equivalence.
\end{thm}

However, while Definition \ref{special} is fairly technical, the condition of speciality is satisfied when, for example, each fibre of the fibration $\bfE \to \bC$ is a model category, and taking a realisation of any simplicial object $X: \Delta^\op \to \cE(c)$ amounts to calculating its homotopy colimit (see \cite{CIS} for the discussion of locally constant functors in this setting). This includes examples like $\DGMod_k$ or any other opfibrations which describe families of model categories with a reasonable notion of geometric realisation leading up to a $\Delta$-structure.

\subsection{Fullness and faithfulness}

The main result of this section which we use to prove Theorem \ref{maintheorem1} is the following proposition: 

\begin{prop}
\label{ff}
Let $F: \cD \to \cC$ be a resolution. Then for any homotopical $\Delta$-opfibration $\cE \to \cC$, the counit transformation
$$
\epsilon: \bF_! \bF^* A \to id_{\Ho \PSect(\cC,\cE)} A
$$
is an isomorphism in $\Ho \PSect(\cC,\cE)$ for any derived section $A$. 
\end{prop}

The proof will be carried out in several steps. Note that for an opfibration $F: \cD \to \cC$ and an object $c \in \cC$, we can take two categories $F / c$, the comma category of $F$ and $c$ (viewed as a functor $[0] \to \cC$), and $\cD(c)$, the fibre of $F$ at $c$. There is a functor which sends $d \in \cD(c)$ to $(d, id_c:F(d) \stackrel = \to c) \in F / c$ and it has a left adjoint given by choosing, for each object $(d, f: F(d) \to c) \in F/c$, an opCartesian morphism $d \to f_! d$ covering $f$. A similar pattern occurs a few times in this section, and this motivates us to introduce the following technical notion:

\begin{opr}
\label{transstructure}
For a category $\cD$, a functor $F: \cD \to \cC$ and an object $c \in \cC$, a \emph{$(F,c)$-transition structure} consists of 
\begin{enumerate}
\item two categories $I, J$ and functors $\mathscr I: I \to \cD$, $\mathscr J: J  \to \cD$,
\item a functor ${\mathrm R}: J \to I$ in $\Cat/\cD$. 
\end{enumerate}
These data are subject to the following conditions:
\begin{itemize}
\item $\rR$ admits a left adjoint ${\mathrm L}$ in $\Cat$,
\item $\mathscr J$ maps $J$ to the fibre $\cD(c)$, so that $F\mathscr J$ factors through $c$.
\end{itemize}
\end{opr}
In the notation of this definition, we sometimes write $(\sI,\sJ,\rR)$ to denote a given $(F,c)$-transition structure.

\begin{exm}
\label{TransitExamples}
The transition structures of importance for us are the following:
\begin{enumerate}
\item For an opfibration $F: \cD \to \cC$ and an object $c \in \cC$, there is a $(F,c)$-transition structure given by  $I = F / c$ and $J = \cD(c)$ outlined just before Definition \ref{transstructure}.  
\item If $F: \cD \to \cC$ is an opfibration and $d \in \cD$, one can have the following $(F, F(d))$-transition structure: $I = \cD / d$ and $J = \cD(F(d))/d$. The right adjoint $\mathrm R$ is the evident inclusion; the left adjoint $\mathrm L$ is given by factoring any morphism $d' \to d$ as 'opCartesian followed by fibrewise' pair of morphisms.
\item Any $(F,c)$ structure $(\mathscr I, \mathscr J, \rR)$ induces a $(F \circ \mathscr I, c)$ structure $(id_I,\rR,\rR)$ with the same right adjoint $\rR$. Thus the first example gives us a $(F_c,c)$-structure where $F_c: \cD/c \to \cC$ is the functor $(d, f: Fd \to c) \mapsto Fd$. For this structure, $I = \cD /c$, $J = \cD(c)$ and $\rR$ acts in the same way as before.
\end{enumerate}
\end{exm}

\begin{rem}
\label{etainduced}
Consider the unit map $\eta(i): i \to \rR \rL i$ for any $i \in I$. Apply $F \circ \mathscr I$ to this map and obtain $\bar \eta(i): F \mathscr I(i) \to c$. For any opfibration $\cE \to \cC$ we then have a well-defined 'restriction' functor
$$
R_c: \mathscr I^* \bF^* \bfE \to \bfE(c) = \cE(c)
$$
where we denote by the same letter $\mathscr I$ the induced functor $\bI \to \bD$. Concretely, this functor sends $(\bi_{[n]}, x \in \cE(F \mathscr I (i_n)))$ to $\bar \eta(i)_! x$, using a (chosen) opCartesian lift $x \to \bar \eta(i)_! x$ in $\cE$ covering $\bar \eta(i): F \mathscr I(i_n) \to c$. 
\end{rem}

\begin{const}
Assume given $c_0 \to ... \to c_n = \bc \in \bC$. Denote by $\bc_!$ the natural transition functor
$$
\bc_!: \bfE(c_0) \cong \cE(c_0) \to \cE(c_n) \cong \bfE(c_n).
$$
Consider also the simplicial comma object (Definition \ref{FovrG}) $\bI \ovr \bR$, where $\bR: \bJ \to \bI$ is the simplicial replacement of $\rR$. Using the diagram (\ref{FovrGdiagram}) and postcomposing with functors to $\bD$ we obtain a new diagram
\begin{equation}
\label{IovrR}
\begin{diagram}[nohug]
      &       & \bI \ovr \bR &       &     \\
      & \ldTo^{pr_{\mathbb I}} & \Leftarrow \dTo \Rightarrow	& \rdTo^{pr_{\mathbb J}} &      \\
\mathbb I & \rTo^{\mathscr I}  & \bD                   & \lTo^{\mathscr J}  &  \mathbb J \\
\end{diagram}
\end{equation}
and we henceforth denote the middle map again by $pr_{\bI \ovr \bR}$.

For $B \in \PSect(\cD, \cE)$ and a given $(F,c_0)$-structure, we get the following diagram
\begin{diagram}[nohug]
      &       & \bI \ovr \bR &       &     \\
      & \ldTo^{pr_{\bI}} & \Leftarrow \dTo \Rightarrow	& \rdTo^{pr_{\bJ}} &      \\
\bI & \rTo_{\bc_! R_{c_0} \mathscr I^* B}  & \bfE(\bc)                   & \lTo_{ \bc_!  \mathscr J^* B}  &  \bJ \\
\end{diagram}
with the middle map $\bc_! pr^*_{\bI \ovr \bR} B$.
Thus we have the span
$$
pr^*_{\bI} \bc_! R_{c_0} \mathscr I^* B \longleftarrow \bc_! pr^*_{\bI \ovr \bR} B \longrightarrow pr^*_{\bJ} \bc_! \mathscr J^* B.
$$
Pushing this forward to the span given by projections,
$$
\Delta^\op \stackrel{\pi_1}{\longleftarrow} \Delta^{\op} \times \Delta^{\op} \stackrel{\pi_2}{\longrightarrow} \Delta^\op,
$$
we obtain a span of bisimplicial objects in $\bfE(c)$:
\begin{equation}
\label{span}
\pi^*_1 \Pi ( \bc_! R_{c_{0}} \mathscr I^*  B) \longleftarrow  \Pi \bc_! pr^*_{\bI \ovr \bR}  B   \longrightarrow \pi^*_2 \Pi (\bc_! \mathscr J^* B).
\end{equation}
Here $\Pi$ is the pushforward functor, simplicial (along $\bI \to \Delta^\op$ and same for $\bJ$)  or bisimplicial (along $\bI \ovr \bR \to \Delta^\op \times \Delta^\op$). We implicitly used the Beck-Chevalley morphisms, such as $\Pi pr^*_{\bI} \to \pi_1^* \Pi$, for pull-backs and pushforwards; they arise from commutative squares like
\begin{diagram}[small]
\bI \ovr \bR & \rTo^{pr_\bI} & \bI \\
\dTo  & 	& \dTo  \\
\Delta^\op \times \Delta^\op & \rTo^{\pi_1} & \Delta^\op \\
\end{diagram} 
\end{const}
by taking associated pull-back functors on functor categories and then replacing some of them by left adjoints.

\begin{rem}
Let us write the terms of the span (\ref{span}) explicitly. For $\bc = c_0 \stackrel{f_1}{\to} ... \stackrel{f_n}{\to} c_n$, we find that
$$
\Pi ( \bc_! R_{c_{0}} \mathscr I^*  B)_m = \coprod_{\bi_{[m]}} (f_n ... f_1 \bar \eta(i_m))_! B(\mathscr I \bi_{[m]}),
$$
where $\bar \eta(i_m)$ is induced from the unit of $\rL \dashv \rR$ (Remark \ref{etainduced}). Next,
$$
\Pi ( \bc_!  \mathscr J^* B)_l = \coprod_{\bj_{[l]}} (f_n ... f_1)_! B(\mathscr J \bj_{[l]}),
$$
and, finally,
$$
\displaystyle (\Pi \bc_! pr^*_{\bI \ovr \bR} B)_{ml} = \coprod_{\bi_{[m]}, \, \, \bj_{[l]}, \, \, \alpha: i_m \to \rR j_0} (f_n ... f_1)_! B( \mathscr I(\bi_{[m]})*^{\mathscr I \alpha} \mathscr J(\bj_{[l]})).
$$
\end{rem}

\begin{prop}
\label{transition}
For $c_0 \stackrel{f_1}{\to} ... \stackrel{f_n}{\to} c_n = \bc  \in \bC$, a $(F:\cD \to \cC, c_0)$-transition structure $(\mathscr I, \mathscr J, \rR)$, and any $B \in \DSect(\cD,F^*\cE)$ there is a natural (in $B$) span of \emph{weak equivalences} in $\bfE(\bc)$
\begin{equation}
\label{TS1}
|\Pi ( \bc_! R_{c_{0}} \mathscr I^*  B)|  \longleftarrow ||\Pi \bc_! pr^*_{\bI \ovr \bR}  B || \longrightarrow |\Pi (\bc_! \mathscr J^* B)|
\end{equation}
which comes from a natural (in $B$) span (\ref{span}) of bisimplicial objects in $\bfE(\bc)$.
\end{prop}
\proof We need to prove that after realisations, both arrows become equivalences. Consider the bisimplicial object
$$
\displaystyle \Pi \pi^*_1 ( \bc_! R_{c_{0}} \mathscr I^*  B)_{ml} = \coprod_{\bi_{[m]}, \, \, \bj_{[l]}, \, \, \alpha: i_m \to \rR j_0} (f_n ... f_1 \bar \eta(i_m))_! B( \mathscr I(\bi_{[m]}))
$$

Our left hand side map in (\ref{span}) passes through this object, as it is equal to the composition
\begin{equation}
\label{TS3}
\Pi \bc_! pr^*_{\bI \ovr \bR} B \to \Pi \pi^*_1 ( \bc_! R_{c_{0}} \mathscr I^*  B) \to \pi^*_1 \Pi ( \bc_! R_{c_{0}} \mathscr I^*  B).
\end{equation}
Writing down the simplicial objects explicitly, we see that the first map in (\ref{TS3}) arises from the action of $B$ on anchor maps and is a termwise weak equivalence of bisimplicial objects because $B$ is a derived section. Realising the second map $\Pi \pi^*_1 ( \bc_! R_{c_{0}} \mathscr I^*  B) \to \pi^*_1 \Pi ( \bc_! R_{c_{0}} \mathscr I^*  B)$ in (\ref{TS3}) along the second simplicial argument, we obtain a map in $\Delta^\op \bfE(\bc)$, whose $m$-th component is 
\begin{equation}
\label{TS4}
|\Pi \pi^*_1 ( \bc_! R_{c_{0}} \mathscr I^*  B)|_m \cong N(i_m \backslash \rR) \otimes  \Pi ( \bc_! R_{c_{0}} \mathscr I^*  B)_m \to  \Pi ( \bc_! R_{c_{0}} \mathscr I^*  B)_m.
\end{equation} 
Observe that because of the adjunction $\rL \dashv \rR$, the category $i_m \backslash \rR=(\rR^{\op} / i_m)^{\op}$ has an initial object (the unit at $i_m$) and is thus contractible, so the map (\ref{TS4}) and thus (\ref{TS3}) and the left-hand side map in (\ref{span}) are all weak equivalences. 

We now have to prove that the right-hand side map $\Pi \bc_! pr^*_{\bI \ovr \bR}  B   \longrightarrow \pi^*_2 \Pi (\bc_! \mathscr J^* B)$ in (\ref{span}) becomes an equivalence after realisations. For each fixed $\bj_{[l]}$, we have a map of simplicial objects, written in components as
\begin{equation}
\label{TS5}
\coprod_{\bi_{[m]}, \, \, \alpha: i_m \to \rR j_0} (f_n ... f_1)_! B( \mathscr I(\bi_{[m]})*^{\mathscr I \alpha} \mathscr J(\bj_{[l]})) \longrightarrow (f_n ... f_1)_! B(\mathscr J(\bj_{[l]})); 
\end{equation} 
because $\rL / j_0$ has a terminal object, Lemma \ref{terminalobject} and Lemma \ref{degeneraciestoequivalences} imply that the map (\ref{TS5}) is a weak equivalence after being realised. We conclude that the morphism 
$$
|\Pi \bc_! pr^*_{\bI \ovr \bR}  B|    \longrightarrow |\pi^*_2 \Pi (\bc_! \mathscr J^* B)|  \cong \Pi (\bc_! \mathscr J^* B)
$$
is an equivalence of simplicial objects in $\bfE(\bc)$, where we took the realisation of bisimplicial objects along the first argument. Proposition \ref{propertiesofhomdelta} then implies that the double realisation
$$
||\Pi \bc_! pr^*_{\bI \ovr \bR}  B||   \longrightarrow ||\pi^*_2 \Pi (\bc_! \mathscr J^* B)| \cong |\Pi (\bc_! \mathscr J^* B)|,
$$
taken in any order, is a weak equivalence.
\endproof

We are now ready to prove Proposition \ref{ff}. Fix $\bc_{[n]} \in \bC$. For $A \in \DSect(\cC,\cE)$ there are functors $A^{aug}_\bc$ and $A_\bc$ (cf. the proof of Lemma \ref{counit}):
\begin{gather*}
A^{aug}_\bc: \bF \ovr c_0 \to \bfE(\bc_{[n]}), \, \, \, \, (\bd_{[m]}, \alpha: Fd_m \to c_0) \mapsto A(\bF(\bd_{[m]})*^\alpha \bc_{[n]}), \\
A_\bc: \bF \ovr c_0 \to \bfE(\bc_{[n]}), \, \, \, \, (\bd_{[m]}, \alpha: Fd_m \to c_0) \mapsto A(\bc_{[n]}). 
\end{gather*}
There is an obvious natural transformation $A^{aug}_\bc \to A_\bc$. Pushing it forward to $\Delta^\op$ and realising gives us a map $A^\bF(\bc) \to N(F /c_0) \otimes A(\bc)$ so that the obvious composition 
$$
A^\bF(\bc) \to N(F /c_0) \otimes A(\bc) \to A(\bc)
$$ 
is the $\bc$-th component of the right-hand map of the counit correspondence (\ref{corrcounit}).

\begin{lemma}
The morphism $N(F /c_0) \otimes A(\bc) \to A(\bc)$ is a weak equivalence.
\end{lemma}
\proof There is an adjunction $F / c_0 \rightleftharpoons \cD(c_0)$ and $\cD(c_0)$ is contractible, thus $F / c_0$ is contractible as well because adjunctions of categories are known to induce homotopy equivalences between the associated nerves \cite{Q}. \endproof

Now recall Example \ref{TransitExamples}(3) where we work over $F/c_0$, with $I = F/c_0$, $J = \cD(c_0)$ and $\rR: \cD(c_0) \to F/c_0$ being the evident functor. Also take the trivial opfibration $\cE(c_n) \times \cC \to \cC$. Both $A_\bc$ and $A_\bc^{aug}$ are then sections over $\bF \ovr c_0$ of the trivial fibration $\bfE(\bc_{[n]}) \times \bF \ovr c_0 \to \bF \ovr c_0$.

\begin{lemma}
The map $A^\bF(\bc) \to N(F /c_0) \otimes A(\bc)$ is a weak equivalence.
\end{lemma}
\proof The obvious natural transformation $A^{aug}_\bc \to A_\bc$, when plugged in the left hand side of the span (\ref{TS1}) for the transition structure of the Example \ref{TransitExamples}(3), gives us the map in question, $A^\bF(\bc) \to N(F /c_0) \otimes A(\bc)$. The right-hand side of span (\ref{TS1}) gives the map
\begin{equation}
\label{TS6}
|\Pi (\rR^*A^{aug}_\bc)| \to N(\cD(c_0)) \otimes A(\bc)
\end{equation}
so by Proposition \ref{transition} we are done if the map (\ref{TS6}) is a weak equivalence. Observe however that 
$$
\Pi (\rR^*A^{aug}_\bc)_m = \coprod_{\bd_{[m]} \in \bD(c_0)} A(\bF(\bd_{[m]}) * \bc) = \coprod_{\bd \in \bD(c_0)} A(id^m_{c_0} * \bc)
$$
with $id_{c_0}^m$ being the degenerate $m$-simplex $c_0 \stackrel{id_{c_0}}{\to} ... \stackrel{id_{c_0}}{\to} c_0$. Because $A$ is a derived section, Lemma \ref{degeneraciestoequivalences} implies that the obvious map $A(id_{c_0}^m * \bc) \to A(\bc)$ is a weak equivalence, so that 
$$
\Pi (\rR^*A^{aug}_\bc)_m \to N(\cD(c_0))_m \otimes A(\bc) = \Pi (\rR^* A_\bc)_m
$$
is a weak equivalence as well.
\endproof

Varying $\bc$, we obtain the proof of Proposition \ref{ff}. With Corollary \ref{ffcondition}, we get that $\bF^*$ is fully faithful on homotopy level, which is exactly the contents of Theorem \ref{maintheorem1}. 

\subsection{Essential surjectivity}

Our second main result, Theorem \ref{maintheorem2}, needs a technical condition of speciality. To state it, we need to define a few auxiliary things. First, take any opfibration $F:\cD \to \cC$. When $F$ is viewed as a functor $\cC \to \Cat$, we can compose it with the endofunctor $\Cat \to \Cat$ which is the simplicial replacement functor. On the level of opfibrations, define the category\footnote{The dependence of the definition of $\bO_\cC (\cD)$ on $F$ is implicit in the notation.} $\bO_\cC (\cD)$ as follows. An object of $\bO_\cC (\cD)$ is an object $c \in \cC$ and $\bd  \in \bD(c)$. A morphism $(c,\bd_{[n]}) \to (c',\bd'_{[m]})$ consists of a map $f: c \to c'$ and an equivalence class of pairs $(\beta,\gamma)$ where
\begin{itemize}
\item $\beta: \bd_{[n]} \Rightarrow \bd^0_{[n]}$ is some natural transformation in $Fun([n],\cD)$ with domain $\bd_{[n]}$ and so that each $\beta_i: d_i \to d^0_i$ is an opCartesian morphism in $\cD$ lying over $f: c \to c'$,
\item $\gamma: \bd^0_{[n]} \to \bd'_{[m]}$ is a morphism in $\bD(c')$,
\item and the equivalence relation is as follows. Two pairs $(\beta^0: \bd_{[n]} \Rightarrow \bd^0_{[n]}, \gamma^0: \bd^0_{[n]} \to \bd'_{[m]})$ and $(\beta^1: \bd_{[n]} \Rightarrow \bd^1_{[n]}, \gamma^1: \bd^1_{[n]} \to \bd'_{[m]})$ are equivalent if, after applying the functor $\pi: \bD(c') \to \Delta^\op$, we have that $\pi \gamma^0 = \pi \gamma^1$.
\end{itemize}
In all, we obtain an opfibration $\bO_\cC (\cD) \to \cC$ whose fibres are $\bD(c)$ and whose transition functors are given by the simplicial replacements of $f_!$, the transition functors of $F: \cD \to \cC$ associated to $f: c \to c'$.

For any opfibration $F: \cD \to \cC$, denote by $F^*F: F^* \cD \to \cD$ the pull-back opfibration of $F$ along $F$. Then from $F^* F$ we obtain the opfibration $\bO_\cD (F^* \cD) \to \cD$ constructed as above, and denote its pull-back along the first element map $h_\cD: \bD \to \cD$ (see Lemma \ref{headandtail}) by $\bO(F^* \cD) \to \bD$. Finally, take the power fibration
$$
(\bF^*\bfE)^{\bO(F^* \cD)} \to \bD.
$$
The $\Delta$-structure, as usual, gives us the lax realisation morphism
\begin{diagram}[small]
(\bF^*\bfE)^{\bO(F^* \cD)}	&		& \rTo^{|-|} 	&			&	\bF^*\bfE \\
	& \rdTo	&				&	\ldTo	&		\\
	&		&	\bD			&			&	\\
\end{diagram}
defined by taking $X \in (\bF^*\bfE)^{\bO(F^* \cD)}$, which is a functor $\bD(F(d_0)) \to \bfE(\bd_{[n]})$ for some $\bd_{[n]} \in \bD$, and realising it (cf. Definition \ref{realisationofafunctor}). There is also, however, the 'evaluation' map
\begin{diagram}[small]
(\bF^*\bfE)^{\bO(F^* \cD)}	&		& \rTo^{ev} 	&			&	\bF^*\bfE \\
	& \rdTo	&				&	\ldTo	&		\\
	&		&	\bD			&			&	\\
\end{diagram}
given by sending the same $X$ to $X(d_0)$, since $d_0 \in \bD(F(d_0))$. The inclusion $X(d_0) \to \coprod_{d \in \bD(F(d_0))} X(d)$ defines a natural transformation $i:ev \Rightarrow |-|$.
\begin{opr}
\label{special}
Given a resolution $F: \cD \to \cC$, a homotopical $\Delta$-opfibration $\cE \to \cC$ is \emph{$F$-special} iff for each $X \in (\bF^*\bfE)^{\bO(F^* \cD)}$, which, when viewed as a functor $\bD(F(d_0)) \to \bfE(\bd_{[n]})$, sends all maps of $\bD(F(d_0))$ to weak equivalences in $\bfE(\bd)$, the $X$-th component of the natural transformation $i$,
$$
i_X: ev(X) \to |X|.
$$ 
is a weak equivalence  
\end{opr}
 
The result of this section is the following. Let $F : \cD \to \cC$ be a resolution.

\begin{prop}
For a $F$-special (Definition \ref{special}) homotopical $\Delta$-opfibration $\cE \to \cC$ and a locally constant $B \in \DSect(\cD,\cE)$, the map
$$
id_{\bD!} id^*_{\bD} B \to \bF^* \bF_! B
$$ 
is a weak equivalence.  
\end{prop}

We will prove that for each $\bd$, the map $id_{\bD!} id^*_{\bD} B(\bd) \to \bF^* \bF_! B(\bd)$ is an equivalence.

\begin{opr}
Given a $(F: \cD \to \cC,c)$-structure $(\mathscr I,\mathscr J,\rR)$ and a $(F': \cD \to \cC',c')$-structure $(\mathscr I',\mathscr J',\rR')$, a \emph{morphism} from the first to the second one consists of  
\begin{itemize}
\item a functor $G: \cC \to \cC'$ in $\cD \backslash \Cat$ with $G(c)=c'$.
\item a commutative square in $\Cat/\cD$
\begin{diagram}[small]
I 			& \lTo^\rR & J \\
\dTo<\lambda &		&	\dTo>\mu \\
I' 			& \lTo^{\rR'} & J'. \\
\end{diagram}
\end{itemize}  
\end{opr}

\begin{exm}
\label{TransitMapExm}
In Example \ref{TransitExamples}, there is a morphism from the second to the first example as soon as $c = F(d)$. In detail: we have a $(F, F(d))$ transition structure $\rL: \cD/d \rightleftharpoons \cD(F(d))/d: \rR$ and $(F,c=F(d))$ transition structure $\rL':F/c \rightleftharpoons \cD(c): \rR'$. In the notation of the definition, $G$ is simply given by $id_\cC$ (this works because $F(d)=c$), $\lambda$ is given by mapping $\alpha:d' \to d$ to $(d',F(\alpha):F(d') \to c)$ and $\mu$ is the evident functor $\cD(c)/d \to \cD(c)$. In this case, even more is true: the square with left adjoints
\begin{diagram}[small]
\cD/d 			& \rTo^\rL & \cD(c)/d \\
\dTo<\lambda &		&	\dTo>\mu \\
F/c 			& \rTo^{\rL'} & \cD(c). \\
\end{diagram}
commutes up to isomorphism.
\end{exm}

\begin{rem}
\label{pullpush}
Given functors $p: A \to B$, $q: A' \to B$ and $r: A' \to A$ such that $p r = q$, for any other functor $X: A \to \cM$ to a cocomplete category, there is a natural map $q_! r^* X \rightarrow p_! X$ where as usual, $r^*$ denotes pull-back and $p_!, q_!$ denote pushforward functors (left adjoint to pull-backs $p^*$ and $q^*$).
\end{rem}

\begin{lemma}
Fix $\bc \in \bC$ and $\bc' \in \bC'$. Let $(\mathscr I,\mathscr J,\rR)$ be a $(F: \cD \to \cC,c_0)$-structure and $(\mathscr I',\mathscr J',\rR')$ be a $(F': \cD \to \cC',c'_0)$-structure. For any morphism $(G: \cC' \to \cC,\lambda,\mu)$ of these transition structures such that $\bG(\bc') = \bc$, a homotopical $\Delta$-opfibration $\cE \to \cC$ and a presection $B: \bD \to \bF^* \bfE$, there is an induced morphism of spans
\begin{equation}
\label{mapoftrans}
\begin{diagram}
\pi^*_1 \Pi ( \bc'_! R_{c'_{0}} \mathscr I'^*  B) & \lTo &  \Pi \bc'_! pr^*_{\mathbb I' \ovr \mathbb R'}  B   & \rTo & \pi^*_2 \Pi (\bc'_! \mathscr J'^* B) \\
\dTo		&		&	\dTo		&		&	\dTo \\
\pi^*_1 \Pi ( \bc_! R_{c_{0}} \mathscr I^*  B) & \lTo &  \Pi \bc_! pr^*_{\bI \ovr \bR}  B   & \rTo & \pi^*_2 \Pi (\bc_! \mathscr J^* B). \\
\end{diagram}
\end{equation}
\end{lemma}
 
\proof The maps exist due to Remark \ref{pullpush}. To get the rightmost map of (\ref{mapoftrans}), apply $\pi_2^*$ to
$$
\Pi (\bc'_! \mathscr J'^* B) \to \Pi (\bc_! \mathscr J^* B)
$$
which we get due to the fact that $\mu^* \bc_! \mathscr J^* B = \bc'_! \mathscr J'^* B$. The middle map of (\ref{mapoftrans}) is obtained in this way as well, and so is the leftmost map (observe that due to the conditions imposed, both restriction functors $R_{c_0}$ and $R_{c'_0}$ agree).

One can then check the commutativity of the squares obtained through a direct computation. For example, observe that the middle map, in components
$$
\underset{\alpha': i'_m \to \rR j'_0}{\coprod_{\bi'_{[m]}, \, \, \bj'_{[l]},}} \bc'_! B( \mathscr I'(\bi'_{[m]})*^{\mathscr I' \alpha'} \mathscr J'(\bj'_{[l]})) \to \underset{\alpha: i_m \to \rR j_0}{\coprod_{\bi_{[m]}, \, \, \bj_{[l]},}} \bc_! B( \mathscr I(\bi_{[m]})*^{\mathscr I \alpha} \mathscr J(\bj_{[l]}))
$$
is induced by the maps of sets indexing the coproducts, given by $(\bi', \bj', \alpha') \mapsto (\lambda(\bi'), \mu(\bj'), \lambda \alpha')$ (cf. the proof of Lemma \ref{corrlemma3}).   
\endproof

\begin{corr}
\label{CorrTS}
Given a map of two transition structures and a \emph{derived section} $B$, the following are equivalent
\begin{enumerate}
\item $|\Pi ( \bc'_! R_{c'_{0}} \mathscr I'^*  B)| \to |\Pi ( \bc_! R_{c_{0}} \mathscr I^*  B)|$ is a weak equivalence,
\item $|\Pi (\bc'_! \mathscr J'^* B)| \to |\Pi (\bc_! \mathscr J^* B)|$ is a weak equivalence.
\end{enumerate}
\end{corr}
\proof Evident. \endproof
 
We now apply that to Example \ref{TransitMapExm}. Observe that for $\bd \in \bD$ with $\bc = \bF(\bd)$, the map
$$
id_{\bD!} id^*_{\bD} B(\bd) \to \bF^* \bF_! B(\bd)
$$ 
exactly corresponds to the first morphism in Corollary \ref{CorrTS}. Writing $\bd$ instead of $\bc'$, observe that the objects in the second map of Corollary \ref{CorrTS} are
$$
\Pi (\bd_! \mathscr J'^* B)_m = \coprod_{\bd'_{m} \in \bD(F(d_0)), \, \, \, d'_m \to d_0} \bF(\bd)_! B(\bd'_{[m]}),
$$
$$
\Pi (\bc_! \mathscr J^* B)_m = \coprod_{\bd'_{m} \in \bD(F(d_0))} \bF(\bd)_! B(\bd'_{[m]}).
$$
The realisation of the first object is equivalent to $\bF(\bd)_! B(d_0)$. It is easy to check that for a $F$-special homotopical $\Delta$-opfibration the functor
$$
\bD(F(d_0)) \to \bfE(\bF(\bd)), \, \, \, \, \bd' \mapsto \bF(\bd)_! B(\bd') 
$$ 
which sends all morphisms to weak equivalences has its realisation equivalent to $\bF(\bd)_! B(d_0)$ and this implies that the map $|\Pi (\bd_! \mathscr J'^* B)| \to |\Pi (\bc_! \mathscr J^* B)|$ is an equivalence.

\begin{lemma}
\label{reflect}
For $F: \cD \to \cC$ a resolution, the functor $\bF^*$ reflects the condition of being a derived section. That is, if $\bF^* A$ is a derived section, then $A$ is one as well.
\end{lemma}
\proof If $\bF^* A$ is a derived section for $A \in \PSect(\cC,\cE)$, then take any anchor map $\bc' \to \bc$ and find an anchor map $\bd' \to \bd$ such that $\bF(\bd' \to \bd) = \bc' \to \bc$ (this is possible due to $F$ being an opfibration with contractible, and hence non-empty, fibres). Then since $\bF^* A(\bd' \to \bd) = A(\bc' \to \bc)$, we get that $A$ is a derived section. \endproof

\begin{corr}[proof of Theorem \ref{maintheorem2}]
 $\bF_!$ sends locally constant sections to derived sections, and
 $$
 \bF_! : \Ho \DSect_{lc}(\cD,\cE) \rightleftarrows \Ho \DSect(\cC,\cE): \bF^*
 $$
 is an equivalence of categories for a special homotopical $\Delta$-fibration $\cE \to \cC$. 
 \end{corr}
 \proof We proved that the unit correspondence gives an isomorphism $id \to \bF^* \bF_!$ of functors on $\Ho \PSect_{lc}(\cD,\cE)$. Using Lemma \ref{reflect}, we see that then $\bF_!$ preserves derived section condition for locally constant sections. This allows us to restrict the unit $id \to \bF^* \bF_!$ to the derived sections. \endproof

$\, $

{\bf Edouard BALZIN}

{\it Physical Address:} National Research University Higher School of Economics, 7 Vavilova st, 117312 Moscow, Russia, and Laboratoire J. A. Dieudonn\'e, University of Nice, Parc Valrose, 06106 Nice, France.

{\it E-mails}: balzin@unice.fr, erbalzin@edu.hse.ru

\end{document}